\documentclass[12pt]{article}

\usepackage{amssymb}
\usepackage{emlines}
\usepackage{luis}
\input{epsf.sty}

\def\stackunder#1#2{\mathrel{\mathop{#2}\limits_{#1}}}
\def\QTR#1#2{{\csname#1\endcsname #2}}

\def\func#1{{\rm #1}}
\def\limfunc#1{\mathop{\rm #1}}%
\def\dfrac#1#2{{\displaystyle {#1 \over #2}}}

\newtheorem{theorem}{Theorem}[section]

\newtheorem{corollary}[theorem]{Corollary}

\newtheorem{proposition}[theorem]{Proposition}
\newtheorem{definition}[theorem]{Definition}

\newtheorem{remark}[theorem]{Remark}

\begin{document}

\thispagestyle{empty}

\vspace*{3cm}
\begin{center}
{\LARGE Differential Geometry on Compound \\ Poisson Space\bigskip }

{\sc Yuri G.~Kondratiev}$^{1,2,3}$

{\sc Jos\'{e} L.~Silva}$^{2,4}$

{\sc Ludwig Streit}$^{2,4}$
\bigskip 

\begin{tabular}{l}
$^{1}$Inst.~Angewandte Math., Bonn Univ., D-53115 Bonn, Germany\\
$^{2}$BiBoS, Bielefeld Univ., D-33615 Bielefeld, Germany \\ 
$^{3}$Inst.~Math., 252601 Kiev, Ukraine\\ 
$^{4}$CCM, Univ.~Madeira, P-9000 Funchal, Portugal \\ 
\multicolumn{1}{r}{CCM$^{@}$-UMa 24/97}
\end{tabular}
\end{center}

\renewcommand{\thefootnote}{}
\footnotetext{
$^{@}$http:/www.uma.pt/ccm/ccm.html}

\newpage
\renewcommand{\thefootnote}{\arabic{footnote}}
\setcounter{page}{1}

\title{
Differential Geometry on Compound \\ Poisson Space}
\author{\textbf{Yuri G.~Kondratiev} \\
Institut f\"ur Angewandte Mathematik, Bonn \\
Universit\"at, D 53115 Bonn, Germany\\
and\\
Forschungszentrum BiBoS, Bielefeld \\
Universit\"at, D 33615 Bielefeld, Germany \\
and\\            
Institute of Mathematics, NASU, 252601 Kiev, Ukraine\\ 
\and
\textbf{Jos\'{e} L.~Silva} \\
Forschungszentrum BiBoS, Bielefeld \\
Universit\"at, D 33615 Bielefeld, Germany \\            
and\\
Centro de Ci\^encias Mathematics, Universidade da Madeira,\\
P 9000 Funchal, Portugal (luis@uma.pt)\\
\and 
\textbf{Ludwig Streit} \\
Forschungszentrum BiBoS, Bielefeld \\
Universit\"at, D 33615 Bielefeld, Germany \\            
and\\
Centro de Ci\^encias Matem\'aticas, Universidade da Madeira,\\
P 9000 Funchal, Portugal
}
\date{}
\maketitle
\renewcommand{\thefootnote}{}
\footnotetext{
Published in {\textit{Methods of Functional Analysis and Topology}} 
{\bf 4}(1), pp. 32--58, 1998.
}

\thispagestyle{empty}

\newpage
\renewcommand{\thefootnote}{\arabic{footnote}}
\setcounter{page}{1}

\newpage
\begin{abstract}
In this paper we carry out analysis and geometry for a class of infinite
dimensional manifolds, namely, compound configuration spaces as a natural
generalization of the work \cite{AKR97}. More precisely a differential
geometry is constructed on the compound configuration space $\Omega _{X}$
over a Riemannian manifold $X.$ This geometry is obtained as a natural
lifting of the Riemannian structure on $X$. In particular, the intrinsic
gradient $\nabla ^{\Omega _{X}}$ divergence \textrm{div}$_{\pi _{\sigma
}^{\tau }}^{\Omega _{X}}$, and Laplace-Beltrami operator $H_{\pi _{\sigma
}^{\tau }}^{\Omega _{X}}=-$\textrm{div}$_{\pi _{\sigma }^{\tau }}^{\Omega
_{X}}\nabla ^{\Omega _{X}}$ are constructed. Therefore the corresponding
Dirichlet forms $\mathcal{E}_{\pi _{\sigma }^{\tau }}^{\Omega _{X}}$ on $%
L^{2}(\Omega _{X},\pi _{\sigma }^{\tau })$ can be defined. Each is shown to
be associated with a diffusion process on $\Omega _{X}$ (so called
equilibrium process) which is nothing but the diffusion process on the
simple configuration space $\Gamma _{X}$ together with corresponding marks,
i.e, $(X_{t}^{\gamma _{\omega }},m_{\omega })$. As another consequence of
our results we obtain a representation of the Lie-algebra of compactly
supported vector fields on $X$ on compound Poisson space. Finally
generalizations to the case when $\pi _{\sigma }^{\tau }$ is replaced by a
marked Poisson measure $\mu _{\sigma \otimes \tau }$ easily follow from
this construction.
\end{abstract}

\newpage
\tableofcontents

\section{Introduction}

Starting with the work of Gelfand et al.~\cite{GGV75}, many researchers
consider representations on compound Poisson space, see also \cite{I96}.
Hence it is natural to ask about geometry and analysis on this space. On the
other hand in statistical physics of continuous systems compound Poisson
measures and their Gibbsian perturbation are used for the description of
many concrete models, see e.g.~\cite{AGL78}.

In constructing analysis and geometry in the space of simple configurations $%
\Gamma _{X}$ over a manifold $X$, i.e., 
\[
\Gamma _{X}=\{\gamma \subset X\,|\,|\gamma \cap K|<\infty \;\mathrm{%
for\;any\;compact\;}K\subset X\}, 
\]
an important tool is the action of the group of diffeomorphism \textrm{Diff}$%
_{0}(X)$ on $X$ which are equal to the identity outside a compact on the
configuration space $\Gamma _{X}$ (similar to shifts on linear spaces), cf.~%
\cite{AKR97}.

In this paper we present a natural extension of the results obtained in \cite
{AKR97} to the case of compound configuration space $\Omega _{X}$ over a
Riemannian manifold $X$, i.e., the space of $\Bbb{R}_{+}$-valued measures on 
$X$ of the form 
\[
\Omega _{X}=\{\omega =\sum_{x\in \gamma _{\omega }}s_{x}\varepsilon _{x}\in 
\mathcal{D}^{\prime }|s_{x}\in \mathrm{supp\,}\tau ,\,\gamma _{\omega }\in
\Gamma _{X}\}, 
\]
where $\tau $ is a finite measure on $\Bbb{R}_{+}$. This geometry is
constructed via a ``lifting procedure'' and is completely determined by the
Riemannian structure on $X$ (cf.~Subsection \ref{3eq58}). In particular, we
obtain the corresponding intrinsic gradient $\nabla ^{\Omega _{X}}$,
divergence \textrm{div}$_{\pi _{\sigma }^{\tau }}^{\Omega _{X}}$, and
Laplace-Beltrami operator $\triangle ^{\Omega _{X}}=$\textrm{div}$_{\pi
_{\sigma }^{\tau }}^{\Omega _{X}}\nabla ^{\Omega _{X}}$. For details we
refer to Section \ref{3eq59} and Section \ref{3eq60}. Here we only mention
that the ``tangent bundle'' $T\Omega _{X}$ of $\Omega _{X}$ is given as
follows 
\[
T_{\omega }\Omega _{X}:=L^{2}(X\rightarrow TX;\omega ),\;\omega \in \Omega
_{X}, 
\]
i.e., the space of sections in the tangent bundle $TX$ of $X$ which are
square-integrable with respect to the Radon measure $\omega $. Since each $%
T_{\omega }\Omega _{X}$ is thus a Hilbert space (endowed with the
corresponding $L^{2}$-inner product $\langle \cdot ,\cdot \rangle
_{T_{\omega }\Omega _{X}}$ coming from the measure $\omega $) $\Omega _{X}$
obtains a Riemannian-type structure which is non-trivial (i.e., varies with $%
\omega $) even when $X=\Bbb{R}^{d}$.

The problem of analysis and geometry on infinite dimensional spaces is
highly connected with the lack of a good notion of ``volume element'' which
is due to the fact that there is no Lebesgue measure on infinite dimensional
linear spaces. The volume element on $X$ is (up to constant multiples) the
unique positive Radon measure $\mu $ on $X$ such that the gradient $\nabla
^{X}$ and the divergence \textrm{div}$^{X}$ become dual operators on $%
L^{2}(X,\mu )$ (w.r.t.~$\langle \cdot ,\cdot \rangle _{TX}$), see e.g.~\cite
{C84}. In Subsection \ref{3eq61} we prove that the probability measure $\pi
_{\sigma }^{\tau }$ on $\Omega _{X}$ for which $\nabla ^{\Omega _{X}}$ and 
\textrm{div}$_{\pi _{\sigma }^{\tau }}^{\Omega _{X}}$ become dual operators
on $L^{2}(\Omega _{X},\pi _{\sigma }^{\tau })$ (w.r.t.~$\langle \cdot ,\cdot
\rangle _{T\Omega _{X}}$) is the right ``volume element'' corresponding to
our differential geometry on $\Omega _{X}$. Of course for completeness
concerning the ``volume element'' one should investigate for which class of
measures the above result is valid but we do not explore this question in
the paper.

Let us stress that the ``test'' functions $\mathcal{F}C_{b}^{\infty }(%
\mathcal{D},\Omega _{X})$ (resp.~``test'' vector fields $V$) we consider as
domains for our gradient $\nabla ^{\Omega _{X}}$ (resp.~\textrm{div}$_{\pi
_{\sigma }^{\tau }}^{\Omega _{X}}$) above are of cylinder type, i.e., $F\in 
\mathcal{F}C_{b}^{\infty }(\mathcal{D},\Omega _{X})$ if and only if 
\[
\omega \mapsto F(\omega )=g_{F}(\langle \omega ,\varphi _{1}\rangle ,\ldots
,\langle \omega ,\varphi _{N}\rangle ) 
\]
for some $N\in \Bbb{N}$, $\varphi _{1},\ldots ,\varphi _{N}\in \mathcal{D}%
:=C_{0}^{\infty }(X)$, $g_{F}\in C_{b}^{\infty }(\Bbb{R}^{N})$ (and $V$
correspondingly, cf.~(\ref{3eq24})). Hence so far the analysis on $\Omega _{X}
$ is basically finite dimensional. However, we can do generic infinite
dimensional analysis on $\Omega _{X}$ by introducing the first order Sobolev
space $H_{0}^{1,2}(\Omega _{X},\pi _{\sigma }^{\tau })$ by closing the
corresponding Dirichlet form 
\[
\mathcal{E}_{\pi _{\sigma }^{\tau }}^{\Omega _{X}}(F,G)=\int_{\Omega
_{X}}\langle (\nabla ^{\Omega _{X}}F)(\omega ),(\nabla ^{\Omega
_{X}}G)(\omega )\rangle _{T_{\omega }\Omega _{X}}d\pi _{\sigma }^{\tau
}(\omega ), 
\]
on $L^{2}(\Omega _{X},\pi _{\sigma }^{\tau })$, i.e., a function $F\in
H_{0}^{1,2}(\Omega _{X},\pi _{\sigma }^{\tau })$ is together with its
gradient $\nabla ^{\Omega _{X}}F$ obtained as a limit in $L^{2}(\Omega ,\pi
_{\sigma }^{\tau })$ of a sequence $F_{n}\in \mathcal{F}C_{b}^{\infty }(%
\mathcal{D},\Omega _{X})$, resp.~$\nabla ^{\Omega _{X}}F_{n}$, $n\in \Bbb{N}$%
. Thus such $F$ really depends on infinitely many points in $X$
(cf.~Subsections \ref{3eq62}, \ref{3eq63}).

Since on $X$ there is a natural diffusion process intrinsically determined
by the geometry, namely distorted Brownian motion on $X$, it is natural to
ask whether the same is true for our geometry on $\Omega _{X}$. In the case
of the space of simple configurations $\Gamma _{X}$ this process is
constructed using the standard theory of Dirichlet forms (cf.~\cite{MR92}),
see recent results in \cite{Y96} and \cite{AKR97}. In the case of compound
configuration space this process has a simple relation with the one
mentioned above. Let us explain this more precisely. We regard every
compound configuration $\omega \in \Omega _{X}$ as depending on two
variables, namely $\omega =(\gamma _{\omega },m_{\omega })$ (or more general
marked configuration) and this allowed us to obtain the following embedding 
\[
L^{2}(\Gamma _{X},\pi _{\sigma })\hookrightarrow L^{2}(\Omega _{X},\pi
_{\sigma }^{\tau }).
\]
As a result we may apply operators acting on $L^{2}(\Gamma _{X},\pi _{\sigma
})$, e.g.~$\nabla ^{\Gamma }$, $\nabla ^{\Gamma *}$ to the space $%
L^{2}(\Omega _{X},\pi _{\sigma }^{\tau })$ acting on part of the variables,
see e.g.~\cite{BK88}. It turns out that the following equality holds 
\[
(\nabla ^{\Omega }F)(\omega )=(\nabla ^{\Gamma }F)((\gamma _{\omega
},m_{\omega })),\;\omega =(\gamma _{\omega },m_{\omega })\in \Omega _{X},
\]
and from this relation it is not hard to obtain relations between the
Dirichlet operators as well as between the correspondings semigroups.
Therefore the process associated to our Dirichlet form is nothing but the
process $X_{t}^{\gamma _{\omega }}$, $t\geq 0$, together with marks, i.e., 
\[
\Xi _{t}=(X_{t}^{\gamma _{\omega }},m_{\omega }),\;t\geq 0,
\]
where $X_{t}^{\gamma _{\omega }}$ is just the equilibrium process on $\Gamma
_{X}$, see \cite[Sect.~6]{AKR97} for details. We describe this procedure in
detail in Section \ref{3eq64}.

We would like to emphasize the contents of Subsection \ref{3eq65}. It is
well-known, see \cite[Sect.~6]{GGV75} that there is a canonical unitary
representation on compound Poisson space, i.e., $L^{2}(\Omega _{X},\pi
_{\sigma }^{\tau })$, of the group of diffeomorphisms \textrm{Diff}$_{0}(X)$%
. On the basis of our results described above, we provide a corresponding
representation of the associated Lie algebra of compactly supported vector
fields. We also exhibit explicit formulas for the corresponding generators.

Thus the contents of Sections \ref{3eq59}, \ref{3eq60}, and \ref{3eq64} have
been described. It remains to add that Section \ref{3eq46} consists of
necessary preliminaries to the furthers sections. Finally in Section \ref
{3eq56} we prove in detail the existence of a marked Poisson measure over the
marked Poisson space $\Omega _{X}^{M}$, where $M$ is a complete separable
metric space with a probability measure. Hence all the results obtained in
this paper extend with trivial changes to marked Poisson space.

Last but not least we would like to mention that most of the results
obtained in this paper extend in a natural way (along the lines of the work 
\cite{AKR97a} and \cite{AGL78}) to the case where compound Poisson measures
are replaced by Gibbs measures of Ruelle type.

Part of the results of this paper were presented in the international
conference ``Analysis on infinite-dimensional Lie algebras and groups''
in Marseille September'97.

\section{Measures on configuration spaces\label{3eq46}}

Let $X$ be a connected, oriented $C^{\infty }$ (non-compact) Riemannian
manifold. For each point $x\in X$, the tangent space to $X$ at $x$ will be
denoted by $T_{x}X$; and the tangent bundle endowed with its natural
differentiable structure will be denoted by $TX=\cup _{x\in X}T_{x}X$. The
Riemannian metric on $X$ associates to each point $x\in X$ an inner product
on $T_{x}X$ which we denote by $\left\langle \cdot ,\cdot \right\rangle
_{T_{x}X}$. The associated norm will be denoted by $|\cdot |_{T_{x}X}$. Let $%
m$ denote the volume element.

$\mathcal{O}(X)$ is defined as the family of all open subsets of $X$ and $%
\mathcal{B}(X)$ denotes the corresponding Borel $\sigma $-algebra. $\mathcal{%
O}_{c}(X)$ and $\mathcal{B}_{c}(X)$ denote the systems of all elements in $%
\mathcal{O}(X)$, $\mathcal{B}(X)$ respectively, which have compact closures.

\subsection{The configuration space over a manifold}

The configuration space $\Gamma _{X}$ over the manifold $X$ is defined as
the set of all locally finite subsets (simple configurations) in $X:$%
\[
\Gamma _{X}=\{\gamma \subset X\,|\,|\gamma \cap K|<\infty \;\mathrm{%
for\;any\;compact\;}K\subset X\}. 
\]
Here $|A|$ denotes the cardinality of a set $A$.

We can identify any $\gamma \in \Gamma _{X}$ with the positive
integer-valued Radon measure 
\[
\sum_{x\in \gamma }\varepsilon _{x}\in \mathcal{M}_{p}(X)\subset \mathcal{M}%
(X), 
\]
where $\sum_{x\in \emptyset }\varepsilon _{x}:=$ zero measure and $\mathcal{M%
}(X)$ $($resp. $\mathcal{M}_{p}(X))$ denotes the set of all positive (resp.
positive integer-valued) Radon measures on $\mathcal{B}\left( X\right) $.
The space $\Gamma _{X}$ can be endowed with the relative topology as a
subset of the space $\mathcal{M}(X)$ with the vague topology, i.e., the
weakest topology on $\Gamma _{X}$ such that all maps 
\[
\Gamma _{X}\ni \gamma \mapsto \left\langle \gamma ,f\right\rangle
:=\int_{X}f(x)d\gamma (x)=\sum_{x\in \gamma }f(x) 
\]
are continuous. Here $f\in C_{0}(X)$ (the set of all real-valued continuous
functions on $X$ with compact support). Let $\mathcal{B}(\Gamma _{X})$
denote the corresponding Borel $\sigma $-algebra. $\mathcal{B}(\Gamma _{X})$
is generated by the sets 
\begin{equation}
C_{\Lambda ,n}:=\{\gamma \in \Gamma _{X}\,|\,|\gamma \cap \Lambda |=n\},
\label{3eq1}
\end{equation}
where $\Lambda \in \mathcal{O}_{c}(X)$, $n\in \Bbb{N}_{0}:=\Bbb{N}\cup \{0\}$%
, see e.g. \cite{GGV75} and \cite{S94}. Note that for any $\Lambda \in 
\mathcal{B}(X)$ and all $n\in \Bbb{N}_{0}$ the set $C_{\Lambda ,n}$ is,
indeed, a Borel subset of $\Gamma _{X}$. Sets of the form (\ref{3eq1}) are
called cylinder sets.

For any $B\subset X$ we introduce a function $N_{B}:\Gamma _{X}\rightarrow 
\Bbb{N}_{0}$ such that 
\[
N_{B}(\gamma )=|\gamma \cap B|,\;\gamma \in \Gamma _{X}. 
\]
Then $\mathcal{B}(\Gamma _{X})$ is the smallest $\sigma $-algebra on $\Gamma
_{X}$ such that all the functions $N_{B}$, $B\in $ $\mathcal{B}_{c}(X)$, are
measurable.

\subsection{Poisson measures\label{3eq12}}

For the construction of a Poisson measure on $\Gamma _{X}$ first we need to
fix an intensity measure $\sigma $ on the underlying manifold $X$. We take a
density $\rho >0$ $m$-a.s.~such that $\rho ^{1/2}\in H_{loc}^{1,2}(X)$ and
put $d\sigma (x)=\rho (x)dm(x)$. Here $H_{loc}^{1,2}(X)$ denotes the local
Sobolev space of order 1 in $L_{loc}^{2}(X,m)$. Then $\sigma $ is a
non-atomic Radon measure on $X$, in particular, $\sigma (\Lambda )<\infty $
for all $\Lambda \in \mathcal{B}_{c}(X)$.

There are different ways to define the Poisson measure $\pi _{\sigma }$ with
intensity $\sigma $ on $\Gamma _{X}$, see e.g. \cite{AKR97} and \cite{GV68}.
Here we characterize $\pi _{\sigma }$ by its Laplace transform.

\begin{definition}
The Laplace transform of $\pi _{\sigma }$ is given for $f\in C_{0}(X)$ by 
\begin{eqnarray}
l_{\pi _{\sigma }}(f) &=&\int_{\Gamma _{X}}\exp (\left\langle \gamma
,f\right\rangle )d\pi _{\sigma }(\gamma )  \nonumber \\
&=&\exp \left( \int_{X}(e^{f(x)}-1)d\sigma (x)\right) .  \label{3eq2}
\end{eqnarray}
\end{definition}

Let us mention that (\ref{3eq2}) defines, via Minlos' theorem, a measure $\pi
_{\sigma }$ on a linear space $F(X)$ of generalized functions on $X$, see
e.g.~\cite{GV68}. An additional analysis shows that the support of the
measure $\pi _{\sigma }$ consists of generalized functions of the form $%
\sum_{x\in \gamma }\varepsilon _{x}$, $\gamma \in \Gamma _{X}$, see e.g.~%
\cite{O87} and \cite{S94}, and then $\pi _{\sigma }$ can be considered as a
measure on $\Gamma _{X}$.

\begin{remark}
By the same argument (\ref{3eq2}) holds for any $\mathcal{B}(X)$-measurable
function $f$ with compact support such that $e^{f}$ is $\sigma $-integrable
on \textrm{supp}$f$. A simple limit-argument then implies that (%
\ref{3eq2}) holds for all $f$ such that $e^{f}-1\in L^{1}(\sigma )$.
\end{remark}

\subsection{Compound Poisson measures\label{3eq22}}

Let $\tau \,$be a non-negative finite measure on $\Bbb{R}_{+}:=]0,\infty [$
having all moments finite and satisfying the analyticity property 
\begin{equation}
\exists C>0:\forall n\in \Bbb{N}_{0}\;\int_{0}^{\infty }s^{n}d\tau
(s)<C^{n}n!.  \label{3eq3}
\end{equation}

We denote $\mathcal{D}=C_{0}^{\infty }(X)$ (the set of $C^{\infty }$%
-functions on $X$ with compact support) equipped with the usual topology,
see e.g. \cite[Chap.~2]{Au82}.

\begin{definition}
A measure $\pi _{\sigma }^{\tau }$ on $\mathcal{D}^{\prime }$ is called a 
\textrm{compound Poisson measure} if its Laplace transform is given
for $\varphi \in \mathcal{D}$ by 
\begin{eqnarray}
l_{\pi _{\sigma }^{\tau }}(\varphi ) &=&\int_{\mathcal{D}^{\prime }}\exp
(\left\langle \omega ,\varphi \right\rangle )d\pi _{\sigma }^{\tau }(\omega )
\nonumber \\
&=&\exp \left( \int_{X}\int_{0}^{\infty }(e^{s\varphi \left( x\right)
}-1)d\tau (s)d\sigma (x)\right) ,  \label{3eq37}
\end{eqnarray}
see e.g.~\cite{GGV75}.
\end{definition}

The measure $\pi _{\sigma }^{\tau }$ has the following properties.

\begin{proposition}
\label{3eq68}

\begin{enumerate}
\item  \label{3eq69}$\pi _{\sigma }^{\tau }$ has an analytic Laplace
transform.

\item  $\pi _{\sigma }^{\tau }$ is supported on $\Omega :=\Omega _{X}$, the
space of compound configurations, i.e., 
\[
\Omega =\{\omega =\sum_{x\in \gamma _{\omega }}s_{x}\varepsilon _{x}\in 
\mathcal{D}^{\prime }|s_{x}\in \func{supp}\tau ,\,\gamma _{\omega }\in
\Gamma _{X}\},
\]
in other words $\pi _{\sigma }^{\tau }(\Omega )=1$.

\item \label{3eq70}If $\func{supp}\tau =\{1\}$, i.e., $d\tau (s)=\varepsilon 
_{1}(ds)$, then $\pi _{\sigma }^{\tau }=\pi _{\sigma }$.
\end{enumerate}
\end{proposition}

For the proof of the above proposition we refer to \cite{KSSU97}, \cite{Ob84}%
. A more general case so called marked Poisson measure is worked out in
Section \ref{3eq56}.

\subsection{The isomorphism between Poisson and compound Poisson spaces\label%
{3eq7}}

Let us define a measure $\hat{\sigma}$ on $(X\times \Bbb{R}_{+},\mathcal{B}%
(X\times \Bbb{R}_{+}))$ as the product measure of the measures $\tau $ and $%
\sigma $, i.e., 
\[
d\hat{\sigma}(\hat{x}):=d\tau (s)d\sigma (x),\;\hat{x}=(x,s)\in X\times \Bbb{%
R}_{+}.
\]
Denote by $\hat{\Gamma}$ the set of the locally finite configurations $\hat{%
\gamma}\subset X\times \Bbb{R}_{+}$ such that 
\[
\hat{\gamma}=\sum_{\hat{x}_{i}\in \hat{\gamma}}\varepsilon _{\hat{x}_{i}},\;%
\hat{x}_{i}=(x_{i},s_{i})\in X\times \Bbb{R}_{+},\,x_{i}\neq x_{j},\;i\neq j
\]
and define the Poisson measure $\pi _{\hat{\sigma}}$ with intensity measure $%
\hat{\sigma}$ on $(\hat{\Gamma},\mathcal{B}(\hat{\Gamma}))\,$via its Laplace
transform 
\begin{eqnarray}
l_{\pi _{\hat{\sigma}}}(\hat{\varphi}) &=&\int_{\hat{\Gamma}}\exp
(\left\langle \hat{\gamma},\hat{\varphi}\right\rangle )d\pi _{\hat{\sigma}}(%
\hat{\gamma})  \nonumber \\
&=&\exp \left( \int_{X\times \Bbb{R}_{+}}(e^{\hat{\varphi}\left( \hat{x}%
\right) }-1)d\hat{\sigma}(\hat{x})\right) ,\,\hat{\varphi}\in \mathcal{D}%
(X\times \Bbb{R}_{+}).  \label{3eq4}
\end{eqnarray}

It follows from (\ref{3eq3}) that the Laplace transform $l_{\pi _{\hat{\sigma 
}}}$ is well defined for $\hat{\varphi }(s,x)=p(s)\varphi (x)$ where $%
p(s)=\sum_{k=0}^{m}p_{k}s^{k}$ $(p_{0}\neq 0)$ is a polynomial and $\varphi
\in \mathcal{D}$ (cf. \cite{LRS95}). Let us put $\hat{\varphi }%
(s,x)=s\varphi (x)$, $\varphi \in \mathcal{D}$ in (\ref{3eq4}). Then by (\ref
{3eq2}) we obtain 
\[
l_{\pi _{\sigma }^{\tau }}(\varphi )=l_{\pi _{\hat{\sigma }}}(s\varphi
),\;\varphi \in \mathcal{D}. 
\]

Then it follows that the compound Poisson measure $\pi _{\sigma }^{\tau }$
is the image of $\pi _{\hat{\sigma }}$ under the transformation $\Sigma :%
\hat{\Gamma }\rightarrow \Sigma \hat{\Gamma }=\Omega \subset \mathcal{D}%
^{\prime }$ given by 
\begin{equation}
\hat{\Gamma }\ni \hat{\gamma }\mapsto (\Sigma \hat{\gamma })(\cdot )=\Sigma
\left( \sum_{\hat{x_{i}}\in \hat{\gamma }}\varepsilon _{\hat{x_{i}}}\right)
(\cdot ):=\sum_{(s_{i},x_{i})\in \hat{\gamma }}s_{i}\varepsilon
_{x_{i}}(\cdot )\in \Omega \subset \mathcal{D}^{\prime },  \label{3eq5}
\end{equation}
i.e., $\forall B\in \mathcal{B}(\mathcal{D}^{\prime })$ 
\[
\pi _{\sigma }^{\tau }(B)=\pi _{\sigma }^{\tau }(B\cap \Omega )=\pi _{\hat{%
\sigma }}(\Sigma ^{-1}(B\cap \Omega )), 
\]
where $\Sigma ^{-1}\Delta $ is the pre-image of the set $\Delta $.

The latter equality may be rewritten in the following form 
\[
\int_{\mathcal{D}^{\prime }}1\!\!1_{B}(\omega )d\pi _{\sigma }^{\tau
}(\omega )=\int_{\Omega }1\!\!1_{B}(\omega )d\pi _{\sigma }^{\tau }(\omega
)=\int_{\hat{\Gamma }}1\!\!1_{B}(\Sigma \hat{\gamma })d\pi _{\hat{\sigma }}(%
\hat{\gamma }), 
\]
which is analogous to the well known change of variable formula for the
Lebesgue integral. Namely, for any $h\in L^{1}(\mathcal{D}^{\prime },\pi
_{\sigma }^{\tau })=L^{1}(\Omega ,\pi _{\sigma }^{\tau })$ the function $%
h\circ \Sigma \in L^{1}(\hat{\Gamma },\pi _{\hat{\sigma }})$ and 
\begin{equation}
\int_{\Omega }h(\omega )d\pi _{\sigma }^{\tau }(\omega )=\int_{\hat{\Gamma }%
}h(\Sigma \hat{\gamma })d\pi _{\hat{\sigma }}(\hat{\gamma }).  \label{3eq6}
\end{equation}

\begin{remark}
It is worth noting that there exists on $\Omega $ an inverse map $\Sigma
^{-1}:\Omega \rightarrow \hat{\Gamma}$. And we obtain that $\pi _{\hat{\sigma%
}}$ on $\hat{\Gamma}$ is the image of $\pi _{\sigma }^{\tau }$ on $\Omega $
under the map $\Sigma ^{-1}$, i.e., $\forall \hat{C}\in \mathcal{B}(\hat{%
\Gamma})$, $\pi _{\hat{\sigma}}(\hat{C})=\pi _{\sigma }^{\tau }(\Sigma \hat{C%
})$ or after rewriting 
\[
\int_{\hat{\Gamma}}1\!\!1_{\hat{C}}(\hat{\gamma})d\pi _{\hat{\sigma}}(\hat{%
\gamma})=\int_{\Omega }1\!\!1_{\Sigma \hat{C}}(\omega )d\pi _{\sigma }^{\tau
}(\omega )=\int_{\Omega }1\!\!1_{\hat{C}}(\Sigma ^{-1}\omega )d\pi _{\sigma
}^{\tau }(\omega ).
\]
As before we easily can write the corresponding change of variables formula,
namely for any $\hat{f}\in L^{1}(\hat{\Gamma},\pi _{\hat{\sigma}})$ the
function $\hat{f}\circ \Sigma ^{-1}\in L^{1}(\Omega ,\pi _{\sigma }^{\tau })$
and 
\begin{equation}
\int_{\hat{\Gamma}}\hat{f}(\hat{\gamma})d\pi _{\hat{\sigma}}(\hat{\gamma}%
)=\int_{\Omega }\hat{f}(\Sigma ^{-1}\omega )d\pi _{\sigma }^{\tau }(\omega ).
\label{3eq57}
\end{equation}
\end{remark}

Now we construct a unitary isomorphism $U_{\Sigma }$ between the Poisson
space $L^{2}(\pi _{\hat{\sigma }}):=L^{2}(\hat{\Gamma },\pi _{\hat{\sigma }%
}) $ and the compound Poisson space $L^{2}(\pi _{\sigma }^{\tau
}):=L^{2}(\Omega ,\pi _{\sigma }^{\tau })$. Namely, 
\[
L^{2}(\Omega ,\pi _{\sigma }^{\tau })\ni h\mapsto U_{\Sigma }h:=h\circ
\Sigma \in L^{2}(\hat{\Gamma },\pi _{\hat{\sigma }}) 
\]
and 
\[
L^{2}(\hat{\Gamma },\pi _{\hat{\sigma }})\ni \hat{f}\mapsto U_{\Sigma }^{-1}%
\hat{f}=\hat{f}\circ \Sigma ^{-1}\in L^{2}(\Omega ,\pi _{\sigma }^{\tau }). 
\]
The isometry of $U_{\Sigma }$ and $U_{\Sigma }^{-1}$ follows from (\ref{3eq6}%
) and (\ref{3eq57}), respectively

As a result we have established the following proposition.

\begin{proposition}
\label{3eq55}The map $U_{\Sigma }$ is a unitary isomorphism between the
Poisson space $L^{2}(\pi _{\hat{\sigma}})$ and the compound Poisson space $%
L^{2}(\pi _{\sigma }^{\tau })$.
\end{proposition}

\subsection{The group of diffeomorphisms and compound Poisson measures\label
{3eq71}}

Let us denote the group of all diffeomorphisms on $X$ by $\func{Diff}(X)$
and by $\func{Diff}_{0}(X)$ the subgroup of all diffeomorphisms $\phi
:X\rightarrow X$ with compact support, i.e., which are equal to the identity
outside of a compact set (depending on $\phi $).

For any $f\in C_{0}(X)$ we have a continuous functional 
\[
\Omega \ni \omega \mapsto \left\langle \omega ,f\right\rangle
=\int_{X}f(x)d\omega (x)=\sum_{x\in \gamma _{\omega }}s_{x}f(x)
\]
and given $\phi \in \func{Diff}_{0}(X)$ we have 
\begin{eqnarray*}
\left\langle \phi ^{*}\omega ,f\right\rangle  &=&\int_{X}f(x)d\omega (\phi
^{-1}(x)) \\
&=&\sum_{x\in \gamma _{\omega }}s_{x}f\circ \phi (x) \\
&=&\left\langle \omega ,f\circ \phi \right\rangle .
\end{eqnarray*}

Any $\phi \in \func{Diff}_{0}(X)$ defines (pointwise) a transformation
of any subset of $X$ and, consequently, the diffeomorphism $\phi $ has the
following ``lifting'' from $X$ to $\Omega :$%
\[
\Omega \ni \omega =\sum_{x\in \gamma _{\omega }}s_{x}\varepsilon _{x}\mapsto
\phi ^{*}\omega =\sum_{x\in \gamma _{\omega }}s_{x}\varepsilon _{\phi
(x)}\in \Omega ,
\]
because for any $f\in C_{0}(X)$%
\begin{eqnarray*}
\int_{X}f(x)d(\phi ^{*}\omega )(x) &=&\int_{X}f(\phi (x))d\omega (x) \\
&=&\sum_{x\in \gamma _{\omega }}s_{x}f(\phi (x)) \\
&=&\int_{X}f(y)\sum_{x\in \gamma _{\omega }}s_{x}\varepsilon _{\phi (x)}(dy).
\end{eqnarray*}

This mapping is obviously measurable and we can define the image $\phi
^{*}\pi _{\sigma }^{\tau }$ of the measure $\pi _{\sigma }^{\tau }$ under $%
\phi $ as usually by $\phi ^{*}\pi _{\sigma }^{\tau }=\pi _{\sigma }^{\tau
}\circ \phi ^{-1}$, i.e., 
\[
(\phi ^{*}\pi _{\sigma }^{\tau })(A)=\pi _{\sigma }^{\tau }(\phi
^{-1}(A)),\;A\in \mathcal{B}(\Omega ). 
\]

The following proposition shows that this transformation is nothing but a
change of the intensity measure $\sigma $, and $\tau $ is preserved.

\begin{proposition}
\label{3eq14}For any $\phi \in \func{Diff}_{0}(X)$ we have 
\[
\phi ^{*}\pi _{\sigma }^{\tau }=\pi _{\phi ^{*}\sigma }^{\tau }.
\]
\end{proposition}

\noindent \textbf{Proof.}\ Due to the characterization of the measures it is
enough to compute the Laplace transform of the measure $\phi ^{*}\pi
_{\sigma }^{\tau }$, to show the property.

Let $f\in C_{0}(X)$ be given. Then the Laplace transform of $\phi ^{*}\pi
_{\sigma }^{\tau }$ is given by 
\begin{eqnarray*}
\int_{\Omega }\exp (\left\langle \omega ,f\right\rangle )d(\phi ^{*}\pi
_{\sigma }^{\tau })(\omega ) &=&\int_{\Omega }\exp (\left\langle \omega
,f\right\rangle )d\pi _{\sigma }^{\tau }(\phi ^{-1}(\omega )) \\
&=&\int_{\Omega }\exp (\left\langle \phi ^{*}\omega ,f\right\rangle )d\pi
_{\sigma }^{\tau }(\omega ) \\
&=&\int_{\Omega }\exp (\left\langle \omega ,f\circ \phi \right\rangle )d\pi
_{\sigma }^{\tau }(\omega ) \\
&=&\exp \left( \int_{X}\int_{0}^{\infty }(e^{sf\circ \phi (x)}-1)d\tau
(s)d\sigma (x)\right)  \\
&=&\exp \left( \int_{X}\int_{0}^{\infty }(e^{sf(x)}-1)d\tau (s)d(\phi
^{*}\sigma )(x)\right)  \\
&=&\int_{\Omega }\exp (\left\langle \omega ,f\right\rangle )d\pi _{\phi
^{*}\sigma }^{\tau }(\omega )
\end{eqnarray*}
which is just the Laplace transform of the measure $\pi _{\phi ^{*}\sigma
}^{\tau }$.\hfill $\blacksquare $

For any $\phi \in \func{Diff}_{0}(X)$ we introduce the Radon-Nikodym
density of $\sigma $ as 
\begin{equation}
\left\{ 
\begin{array}{ll}
p_{\phi }^{\sigma }(x) & :={{\dfrac{d(\phi ^{*}\sigma )}{d\sigma }}}(x)={{%
\dfrac{\rho (\phi ^{-1}(x))}{\rho (x)}}}{{\dfrac{dm(\phi ^{-1}(x))}{dm(x)}}}=%
{{\dfrac{\rho (\phi ^{-1}(x))}{\rho (x)}}}J_{m}^{\phi }(x), \\ 
&  \\ 
& \mathrm{\;\;\;\;\;\;\;\;\;if\;}x\in \{0<\rho <\infty \}\cap \{0<\rho \circ
\phi ^{-1}<\infty \}; \\ 
&  \\ 
p_{\phi }^{\sigma }(x) & :=1,\mathrm{\;otherwise,}
\end{array}
,\right.  \label{3eq23}
\end{equation}
where $J_{m}^{\phi }$ is the Jacobian determinant of $\phi $ (with respect
to the Riemannian volume $m$), see e.g. \cite{Boo75}. Note that $p_{\phi
}^{\sigma }(x)\equiv 1$ outside a compact.

The next proposition is a consequence of the Proposition \ref{3eq55},
Skorokhod's theorem on absolute continuity of Poisson measures, see e.g. 
\cite{Sk57}, \cite{T90} and also \cite{S94}. It shows that $\pi _{\sigma
}^{\tau }$ is quasi-invariant with respect to the group $\func{Diff}%
_{0}(X) $.

\begin{proposition}
\label{3eq17}The compound Poisson measure $\pi _{\sigma }^{\tau }$ is
quasi-invariant with respect to the group $\func{Diff}_{0}(X)$ and for any $%
\phi \in \func{Diff}_{0}(X)$ we have $p_{\phi }^{\pi _{\sigma }^{\tau
}}=p_{\phi }^{\pi _{\lambda _{\tau }\sigma }}$, where $\lambda _{\tau }=\tau
(\Bbb{R}_{+})$, i.e., 
\[
p_{\phi }^{\pi _{\sigma }^{\tau }}(\omega )=\frac{d(\phi ^{*}\pi _{\sigma
}^{\tau })}{d\pi _{\sigma }^{\tau }}(\omega )=\prod_{x\in \gamma _{\omega
}}p_{\phi }^{\sigma }(x)\exp \left( \lambda _{\tau }\int_{X}(1-p_{\phi
}^{\sigma }(x))d\sigma (x)\right) .
\]
\end{proposition}

\noindent \textbf{Proof.}\ Given $\phi \in \func{Diff}_{0}(X)$ then $%
\hat{\phi}:=\phi \otimes \func{id}\in \func{Diff}%
(X\times \Bbb{R}_{+})$. Hence having in mind the isomorphism described in
Subsection \ref{3eq7} the Radon-Nikodym density of $\pi _{\sigma }^{\tau }$
with respect to the group $\func{Diff}_{0}(X)$ is given by 
\begin{eqnarray*}
&&p_{\phi }^{\pi _{\sigma }^{\tau }}(\omega )=U_{\Sigma }^{-1}p_{\phi
\otimes \func{id}}^{\pi _{\hat{\sigma}}}(\omega ) \\
&=&\prod_{\hat{x}\in \hat{\gamma}_{\omega }}\frac{d\hat{\sigma}\circ (\phi
\otimes \func{id})^{-1}}{d\hat{\sigma}}(\hat{x})\exp \left( \int_{X\times 
\Bbb{R}_{+}}\left( 1-\frac{d\hat{\sigma}\circ (\phi \otimes \func{id})^{-1}%
}{d\hat{\sigma}}(\hat{x})\right) d\hat{\sigma}(\hat{x})\right)  \\
&=&\prod_{\hat{x}\in \hat{\gamma}_{\omega }}p_{\phi }^{\sigma }(x)\exp
\left( \lambda _{\tau }\int_{X}(1-p_{\phi }^{\sigma }(x))d\sigma (x)\right) 
\\
&=&\prod_{x\in \gamma _{\omega }}p_{\phi }^{\sigma }(x)\exp \left( \lambda
_{\tau }\int_{X}(1-p_{\phi }^{\sigma }(x))d\sigma (x)\right)  \\
&=&p_{\phi }^{\pi _{\lambda _{\tau }\sigma }}(\gamma _{\omega }),
\end{eqnarray*}
where we have used \cite[Proposition 2.2]{AKR97}.\hfill $\blacksquare $

\section{Differential geometry on compound Poisson space\label{3eq59}}

The underlying differentiable structure on $X$ has a natural lifting to the
configuration space $\Omega $. As a result there appear in $\Omega $ objects
such as the gradient, the tangent space etc. Below we describe the
corresponding constructions in details.

\subsection{The tangent bundle of $\Omega $\label{3eq31}}

Let $V(X)$ be the set of all $C^{\infty }$-vector fields on $X$ (i.e.,
smooth sections of $TX$). We will use a subset $V_{0}(X)\subset V(X)$
consisting of all vector fields with compact support. $V_{0}(X)$ can be
considered as an infinite dimensional Lie algebra which corresponds to the
group \textrm{Diff}$_{0}(X)$ in the following sense: for any $v\in V_{0}(X)$
we can construct the flow of this vector field as a collection of mappings $%
\phi _{t}^{v}:X\rightarrow X$, $t\in \Bbb{R}$ obtained by integrating the
vector field.

More precisely, for any $x\in X$ the curve 
\[
\Bbb{R}\ni t\longmapsto \phi _{t}^{v}(x)\in X 
\]
is defined as the solution to the following Cauchy problem 
\[
\left\{ 
\begin{array}{l}
\dfrac{d}{dt}\phi _{t}^{v}(x)=v(\phi _{t}^{v}(x)) \\ 
\\ 
\phi _{0}^{v}(x)=x
\end{array}
\right. . 
\]
That no explosion is possible and $\phi _{t}^{v}$ is well-defined for each $%
t\in \Bbb{R}$, is a consequence of $v\in V_{0}(X)$ (the latter implies that $%
v$ is a complete vector field). The mappings $\{\phi _{t}^{v},t\in \Bbb{R}\}$
form a one-parameter subgroup of diffeomorphisms in the group \textrm{Diff}$%
_{0}(X)$ (see e.g. \cite{Boo75}), that is, 
\[
\begin{array}{l}
\mathrm{1)\;}\forall t\in \Bbb{R}\;\phi _{t}^{v}\in \mathrm{Diff}_{0}(X) \\ 
\\ 
\mathrm{2)\;}\forall t,s\in \Bbb{R}\;\phi _{t}^{v}\circ \phi _{s}^{v}=\phi
_{t+s}^{v}.
\end{array}
\]
Let us fix $v\in V_{0}(X)$. Having the group $\phi _{t}^{v}$, $t\in \Bbb{R}$%
, we can consider for any $\omega \in \Omega $ the curve 
\[
\Bbb{R}\ni t\longmapsto \phi _{t}^{v}(w)\in \Omega . 
\]

\begin{definition}
\label{3eq9}For a function $F:\Omega \rightarrow \Bbb{R}$ we define the 
\textbf{\emph{directional derivative}} along the vector field $v\in V_{0}(X)$
as 
\[
(\nabla _{v}^{\Omega }F)(\omega ):=\frac{d}{dt}F(\phi _{t}^{v*}\omega
)|_{t=0}, 
\]
provided the right hand side exists.
\end{definition}

We note that $\nabla _{v}^{\Omega }F$ is closely related to the concept of
the Lie derivative corresponding to a special class of vector fields on $%
\Omega $, see below.

Let us introduce a special class of smooth functions on $\Omega $ which play
an important role in our considerations below. We introduce $\mathcal{F}%
C_{b}^{\infty }(\mathcal{D},\Omega )$ as the set of all functions $F:\Omega
\rightarrow \Bbb{R}$ of the form 
\begin{equation}
F(\omega )=g_{F}(\left\langle \omega ,\varphi _{1}\right\rangle ,\ldots
,\left\langle \omega ,\varphi _{N}\right\rangle ),\;\omega \in \Omega ,
\label{3eq8}
\end{equation}
where (generating directions) $\varphi _{1},\ldots ,\varphi _{N}\in \mathcal{%
D}$ and $g_{F}(s_{1},\ldots ,s_{N})$ (generating function for $F$) is from $%
C_{b}^{\infty }(\Bbb{R}^{N})$.

For any $F\in \mathcal{F}C_{b}^{\infty }(\mathcal{D},\Omega )$ of the form (%
\ref{3eq8}) and given $v\in V_{0}(X)$ we have 
\begin{eqnarray*}
F(\phi _{t}^{v*}\omega ) &=&g_{F}(\left\langle \phi _{t}^{v*}\omega ,\varphi
_{1}\right\rangle ,\ldots ,\left\langle \phi _{t}^{v*}\omega ,\varphi
_{N}\right\rangle ) \\
&=&g_{F}(\left\langle \omega ,\varphi _{1}\circ \phi _{t}^{v}\right\rangle
,\ldots ,\left\langle \omega ,\varphi _{N}\circ \phi _{t}^{v}\right\rangle )
\end{eqnarray*}
and, therefore, an application of Definition \ref{3eq9} gives 
\begin{equation}
(\nabla _{v}^{\Omega }F)(\omega )=\sum_{i=1}^{N}\frac{\partial g_{F}}{%
\partial s_{i}}(\left\langle \omega ,\varphi _{1}\right\rangle ,\ldots
,\left\langle \omega ,\varphi _{N}\right\rangle )\left\langle \omega ,\nabla
_{v}^{X}\varphi _{i}\right\rangle ,  \label{3eq10}
\end{equation}
where $\nabla _{v}^{X}\varphi $ is the directional (or Lie) derivative of $%
\varphi :X\rightarrow \Bbb{R}$ along the vector field $v\in V_{0}(X)$, i.e., 
\[
(\nabla _{v}^{X}\varphi )(x)=\left\langle \nabla ^{X}\varphi
(x),v(x)\right\rangle _{T_{x}X}, 
\]
where $\nabla ^{X}$ denotes the gradient on $X$.

The expression of $\nabla _{v}^{\Omega }$ on smooth cylinder functions given
by (\ref{3eq10}) motivates the following definition.

\begin{definition}
\label{3eq28}We introduce the tangent space $T_{\omega }\Omega $ to the
configuration space $\Omega $ at the point $\omega \in \Omega $ as the
Hilbert space of measurable $\omega $ - square - integrable sections
(measurable vector fields) $V_{\omega }:X\rightarrow TX$ with the scalar
product 
\begin{equation}
\left\langle V_{\omega }^{1},V_{\omega }^{2}\right\rangle _{T_{\omega
}\Omega }=\int_{X}\left\langle V_{\omega }^{1}(x),V_{\omega
}^{2}(x)\right\rangle _{T_{x}X}d\omega (x)  \label{3eq29}
\end{equation}
$V_{\omega }^{1},V_{\omega }^{2}\in T_{\omega }\Omega $. The corresponding
tangent bundle is 
\[
T\Omega =\bigcup_{\omega \in \Omega }T_{\omega }\Omega . 
\]
\end{definition}

Let us stress that any $v\in V_{0}(X)$ can be considered as a ``constant''
vector field on $\Omega $ such that 
\[
\Omega \ni \omega \longmapsto V_{\omega }(\cdot )=v(\cdot )\in T_{\omega
}\Omega , 
\]
\[
\left\langle v,v\right\rangle _{T_{\omega }\Omega
}=\int_{X}|v(x)|_{T_{x}X}^{2}d\omega (x)<\infty . 
\]

Usually in Riemannian geometry, having the directional derivative and a
Hilbert space as the tangent space we can introduce the gradient.

\begin{definition}
We define the intrinsic gradient of a function $F:\Omega \rightarrow \Bbb{R}$
as the mapping 
\[
\Omega \ni \omega \longmapsto (\nabla ^{\Omega }F)(\omega )\in T_{\omega
}\Omega 
\]
such that for any $v\in V_{0}(X)$%
\begin{equation}
(\nabla _{v}^{\Omega }F)(\omega )=\left\langle (\nabla ^{\Omega }F)(\omega
),v\right\rangle _{T_{\omega }\Omega }.  \label{3eq21}
\end{equation}
\end{definition}

Note that (\ref{3eq21}), in particular, implies that $\nabla _{v}^{\Omega }F$
is the directional derivative along the ``constant'' vector field $v$ on $%
\Omega $. Furthermore, by (\ref{3eq10}) for any $F\in \mathcal{F}%
C_{b}^{\infty }(\mathcal{D},\Omega )$ of the form (\ref{3eq8}) gives 
\begin{equation}
(\nabla ^{\Omega }F)(\omega ;x)=\sum_{i=1}^{N}\frac{\partial g_{F}}{\partial
s_{i}}(\left\langle \omega ,\varphi _{1}\right\rangle ,\ldots ,\left\langle
\omega ,\varphi _{N}\right\rangle )\nabla ^{X}\varphi _{i}(x),\;\omega \in
\Omega ,x\in X.  \label{3eq11}
\end{equation}

\subsection{Integration by parts and divergence on compound Poisson space%
\label{3eq61}}

Let the configuration space $\Omega $ be equipped with the compound Poisson
measure $\pi _{\sigma }^{\tau }$ (cf. Subsection \ref{3eq22}). The set $%
\mathcal{F}C_{b}^{\infty }(\mathcal{D},\Omega )$ is a dense subset in the
space $L^{2}(\Omega ,\mathcal{B}(\Omega ),\pi _{\sigma }^{\tau })=:L^{2}(\pi
_{\sigma }^{\tau })$. For any vector field $v\in V_{0}(X)$ we have a
differential operator in $L^{2}(\pi _{\sigma }^{\tau })$ on the domain $%
\mathcal{F}C_{b}^{\infty }(\mathcal{D},\Omega )$ given by 
\[
\mathcal{F}C_{b}^{\infty }(\mathcal{D},\Omega )\ni F\longmapsto \nabla
_{v}^{\Omega }F\in L^{2}(\pi _{\sigma }^{\tau }). 
\]
Our aim now is to compute the adjoint operator $\nabla _{v}^{\Omega *}$ in $%
L^{2}(\pi _{\sigma }^{\tau })$. It corresponds, of course, to an integration
by parts formula with respect to the measure $\pi _{\sigma }^{\tau }$.

To this end we recall first of all the integration by parts formula for the
measure $\sigma $. The logarithmic derivative of $\sigma $ is given by the
vector field 
\[
X\ni x\longmapsto \beta ^{\sigma }(x):=\frac{\nabla ^{X}\rho (x)}{\rho (x)}%
\in T_{x}X. 
\]
(where as usual $\beta ^{\sigma }:=0$ on $\{\rho =0\}$). For all $\varphi
_{1},\varphi _{2}\in \mathcal{D}$ we have 
\begin{eqnarray}
&&\int_{X}(\nabla _{v}^{X}\varphi _{1})(x)\varphi _{2}(x)d\sigma (x) 
\nonumber \\
&=&-\int_{X}\varphi _{1}(x)(\nabla _{v}^{X}\varphi _{2})(x)d\sigma
(x)-\int_{X}\varphi _{1}(x)\varphi _{2}(x)\beta _{v}^{\sigma }(x)d\sigma (x),
\label{3eq13}
\end{eqnarray}
where 
\begin{equation}
\beta _{v}^{\sigma }(x):=\left\langle \beta ^{\sigma }(x),v(x)\right\rangle
_{T_{x}X}+\mathrm{div}^{X}v(x)  \label{3eq18}
\end{equation}
is the so-called logarithmic derivative of the measure $\sigma $ along the
vector field $v$ and \textrm{div}$^{X}:=$\textrm{div}$_{m}^{X}$ is the
divergence on $X$ with respect to $m$. Analogously, we define \textrm{div}$%
_{\sigma }^{X}$ as the divergence on $X$ with respect to $\sigma $, i.e., 
\textrm{div}$_{\sigma }^{X}$ is the dual operator on $L^{2}(X,\sigma
)=:L^{2}(\sigma )$ of $\nabla ^{X}$. Then on the one hand we can rewrite (%
\ref{3eq13}) as an operator equality on the domain $\mathcal{D}\subset
L^{2}(\sigma ):$%
\[
\nabla _{v}^{X*}=-\nabla _{v}^{X}-\beta _{v}^{\sigma }, 
\]
where the adjoint operator is considered with respect to $L^{2}(\sigma )$.
Note that, obviously, $\beta _{v}^{\sigma }\in L^{2}(\sigma )$ for all $v\in
V_{0}(X)$. On the other hand we have 
\begin{equation}
\mathrm{div}_{\sigma }^{X}=\beta ^{\sigma }.  \label{3eq26}
\end{equation}

Having the logarithmic derivative $\beta _{v}^{\sigma }$ we introduce an
analogous object for the compound Poisson measure.

\begin{definition}
For any $v\in V_{0}(X)$ we define the logarithmic derivative of the compound
Poisson measure $\pi _{\sigma }^{\tau }$ along $v$ as the following function
on $\Omega :$%
\begin{equation}
\Omega \ni \omega \mapsto B_{v}^{\pi _{\sigma }^{\tau }}(\omega
):=\left\langle \gamma _{\omega },\beta _{v}^{\sigma }\right\rangle
=\int_{X}[\left\langle \beta ^{\sigma }(x),v(x)\right\rangle _{T_{x}X}+%
\mathrm{div}^{X}v(x)]d\gamma _{\omega }(x).  \label{3eq25}
\end{equation}
\end{definition}

A motivation for this definition is given by the following integration by
parts formula.

\begin{theorem}
\label{3eq42}For all $F,G\in \mathcal{F}C_{b}^{\infty }(\mathcal{D},\Omega )$
and any $v\in V_{0}(X)$ we have 
\begin{eqnarray}
&&\int_{\Omega }(\nabla _{v}^{\Omega }F)(\omega )G(\omega )d\pi _{\sigma
}^{\tau }(\omega )  \nonumber \\
&=&-\int_{\Omega }F(\omega )(\nabla _{v}^{\Omega }G)(\omega )d\pi _{\sigma
}^{\tau }(\omega )-\int_{\Omega }F(\omega )G(\omega )B_{v}^{\pi _{\sigma
}^{\tau }}(\omega )d\pi _{\sigma }^{\tau }(\omega ),  \label{3eq16}
\end{eqnarray}
or 
\begin{equation}
\nabla _{v}^{\Omega *}=-\nabla _{v}^{\Omega }-B_{v}^{\pi _{\sigma }^{\tau }}
\label{3eq20}
\end{equation}
as an operator equality on the domain $\mathcal{F}C_{b}^{\infty }(\mathcal{D}%
,\Omega )$ in $L^{2}(\pi _{\sigma }^{\tau })$.
\end{theorem}

\noindent \textbf{Proof.}\ Due to Proposition \ref{3eq14} we have that 
\[
\int_{\Omega }F(\phi _{t}^{v}(\omega ))G(\omega )d\pi _{\sigma }^{\tau
}(\omega )=\int_{\Omega }F(\omega )G(\phi _{-t}^{v}\omega )d\pi _{\phi
_{t}^{v*}\sigma }^{\tau }(\omega ). 
\]

Differentiating this equation with respect to $t$ and interchanging $\frac{d%
}{dt}$ with the integrals, by Definition \ref{3eq9} the left hand side
becomes (\ref{3eq16}). To see that the right hand side also coincides with (%
\ref{3eq16}) we note that 
\[
\frac{d}{dt}G(\phi _{-t}^{v}(\omega ))|_{t=0}=-(\nabla _{v}^{\Omega
}G)(\omega )
\]
and (by Proposition \ref{3eq17}) 
\begin{eqnarray*}
&&\frac{d}{dt}\left. \left[ \frac{d\pi _{\phi _{t}^{v*}\sigma }^{\tau }}{%
d\pi _{\sigma }^{\tau }}(\omega )\right] \right| _{t=0} \\
&=&\left. \frac{d}{dt}\left[ \prod_{x\in \gamma _{\omega }}\frac{\rho (\phi
_{-t}^{v}(x))}{\rho (x)}J_{m}^{\phi _{t}^{v}}(x)\right] \right| _{t=0} \\
&&+\frac{d}{dt}\left. \left[ \exp \left\{ \lambda _{\tau }\int_{X}\left( 1-%
\frac{\rho (\phi _{-t}^{v}(x))}{\rho (x)}J_{m}^{\phi _{t}^{v}}(x)\right)
d\sigma (x)\right\} \right] \right| _{t=0}.
\end{eqnarray*}

Using (\ref{3eq18}) and the formula $\frac{d}{dt}[J_{m}^{\phi _{t}^{v}}(x)%
]|_{t=0}=-\mathrm{div}^{X}v(x)$, the latter expressions becomes equal to 
\begin{eqnarray*}
&-&\sum_{x\in \gamma _{\omega }}[\left\langle \beta ^{\sigma
}(x),v(x)\right\rangle _{T_{x}X}+\mathrm{div}^{X}v(x)] \\
&+&\lambda _{\tau }\int_{X}[\left\langle \beta ^{\sigma
}(x),v(x)\right\rangle _{T_{x}X}+\mathrm{div}^{X}v(x)]d\sigma (x) \\
&=&-\sum_{x\in \gamma _{\omega }}\beta _{v}^{\sigma }(x)+\lambda _{\tau
}\int_{X}\beta _{v}^{\sigma }(x)d\sigma (x)=-B_{v}^{\pi _{\sigma }^{\tau
}}(\omega ),
\end{eqnarray*}
where we have used the equality 
\[
\int_{X}\beta _{v}^{\sigma }(x)d\sigma (x)=-\int_{X}(\nabla
_{v}^{X*}1)(x)d\sigma (x)=0. 
\]
This completes the proof.\hfill $\blacksquare $

\begin{definition}
\label{3eq19}For a vector field 
\[
V:\Omega \ni \omega \longmapsto V_{\omega }\in T_{\omega }\Omega 
\]
the intrinsic divergence $\mathrm{div}_{\pi _{\sigma }^{\tau }}^{\Omega }V$
is defined via the duality relation 
\begin{equation}
\int_{\Omega }\left\langle V_{\omega },(\nabla ^{\Omega }F)(\omega
)\right\rangle _{T_{\omega }\Omega }d\pi _{\sigma }^{\tau }(\omega
)=-\int_{\Omega }F(\omega )(\mathrm{div}_{\pi _{\sigma }^{\tau }}^{\Omega
}V)(\omega )d\pi _{\sigma }^{\tau }(\omega )  \label{3eq44}
\end{equation}
for all $F\in \mathcal{F}C_{b}^{\infty }(\mathcal{D},\Omega )$, provided it
exists $($i.e., provided 
\[
F\longmapsto \int_{\Omega }\left\langle V_{\omega },(\nabla ^{\Omega
}F)(\omega )\right\rangle _{T_{\omega }\Omega }d\pi _{\sigma }^{\tau
}(\omega )
\]
is continuous on $L^{2}(\pi _{\sigma }^{\tau }))$.
\end{definition}

The existence of the divergence, of course, requires some smoothness of the
vector field. A class of smooth vector fields on $\Omega $ for which the
divergence can be computed in an explicit form is described in the following
proposition.

\begin{proposition}
For any vector field 
\begin{equation}
V_{\omega }(x)=\sum_{j=1}^{N}G_{j}(\omega )v_{j}(x),\;\omega \in \Omega
,\;x\in X,  \label{3eq24}
\end{equation}
with $G_{j}\in \mathcal{F}C_{b}^{\infty }(\mathcal{D},\Omega )$, $v_{j}\in
V_{0}(X)$, $j=1,\ldots ,N$, we have 
\begin{eqnarray}
(\mathrm{div}_{\pi _{\sigma }^{\tau }}^{\Omega }V)(\omega )
&=&\sum_{j=1}^{N}(\nabla _{v_{j}}^{\Omega }G_{j})(\omega
)+\sum_{j=1}^{N}B_{v_{j}}^{\pi _{\sigma }^{\tau }}(\omega )G_{j}(\omega ) 
\nonumber \\
&=&\sum_{j=1}^{N}\left\langle (\nabla ^{\Omega }G_{j})(\omega
),v_{j}\right\rangle _{T_{\omega }\Omega }+\sum_{j=1}^{N}\left\langle \gamma
,\beta _{v_{j}}^{\sigma }\right\rangle G_{j}(\omega ).  \label{3eq43}
\end{eqnarray}
\end{proposition}

\noindent \textbf{Proof.}\ Due to the linearity of $\nabla ^{\Omega }$ it is
sufficient to consider the case $N=1$, i.e., $V_{\omega }(x)=G(\omega )v(x)$%
. Then for all $F\in \mathcal{F}C_{b}^{\infty }(\mathcal{D},\Omega )$ and
the Definition \ref{3eq19} we have 
\begin{eqnarray*}
\int_{\Omega }(\mathrm{div}_{\pi _{\sigma }^{\tau }}^{\Omega }V)(\omega
)F(\omega )d\pi _{\sigma }^{\tau }(\omega ) &=&-\int_{\Omega }\left\langle
V_{\omega },(\nabla ^{\Omega }F)(\omega )\right\rangle _{T_{\omega }\Omega
}d\pi _{\sigma }^{\tau }(\omega ) \\
&=&-\int_{\Omega }G(\omega )\left\langle v,(\nabla ^{\Omega }F)(\omega
)\right\rangle _{T_{\omega }\Omega }d\pi _{\sigma }^{\tau }(\omega ) \\
&=&-\int_{\Omega }G(\omega )(\nabla _{v}^{\Omega }F)(\omega )d\pi _{\sigma
}^{\tau }(\omega ) \\
&=&-\int_{\Omega }(\nabla _{v}^{\Omega *}G)(\omega )F(\omega )d\pi _{\sigma
}^{\tau }(\omega ) \\
&=&\int_{\Omega }(\nabla _{v}^{\Omega }G)(\omega )F(\omega )d\pi _{\sigma
}^{\tau }(\omega ) \\
&&+\int_{\Omega }B_{v}^{\pi _{\sigma }^{\tau }}(\omega )G(\omega )F(\omega
)d\pi _{\sigma }^{\tau }(\omega ),
\end{eqnarray*}
where we have used (\ref{3eq20}). Hence 
\begin{eqnarray*}
(\mathrm{div}_{\pi _{\sigma }^{\tau }}^{\Omega }V)(\omega ) &=&(\nabla
_{v}^{\Omega }G)(\omega )+B_{v}^{\pi _{\sigma }^{\tau }}(\omega )G(\omega )
\\
&=&\left\langle (\nabla ^{\Omega }G)(\omega ),v\right\rangle _{T_{\omega
}\Omega }+\left\langle \gamma ,\beta _{v}^{\sigma }\right\rangle G(\omega ).
\end{eqnarray*}
\hfill $\blacksquare $

In the next subsection we give an equivalent description via a ``lifting
rule'' of the above differential objects on $\Omega $. Before we would like
to consider the special case when \textrm{supp}$\tau =\{1\}$, i.e., $d\tau
(s)=\varepsilon _{1}(ds)$, cf.~Proposition \ref{3eq68}-\ref{3eq70}. A
detailed description can be found in \cite{AKR97}. In this case our
intrinsic gradient $\nabla ^{\Omega }$ is nothing but the intrinsic gradient 
$\nabla ^{\Gamma }$ on $L^{2}(\Gamma _{X},\pi _{\sigma })$, i.e., for any $%
F\in \mathcal{F}C_{b}^{\infty }(\mathcal{D},\Gamma )$ of the form 
\[
F(\gamma )=g_{F}(\left\langle \gamma ,\varphi _{1}\right\rangle ,\ldots
,\left\langle \gamma ,\varphi _{N}\right\rangle ),\;\gamma \in \Gamma _{X}, 
\]
where $\varphi _{1},\ldots ,\varphi _{N}\in \mathcal{D}$ and $g_{F}\in
C_{b}^{\infty }(\Bbb{R}^{N})$, we have 
\[
(\nabla ^{\Gamma }F)(\gamma ;x)=\sum_{j=1}^{N}\frac{\partial g_{F}}{\partial
s_{i}}(\left\langle \gamma ,\varphi _{1}\right\rangle ,\ldots ,\left\langle
\gamma ,\varphi _{N}\right\rangle )\nabla ^{X}\varphi _{i}(x),\;\gamma \in
\Gamma _{X},x\in X. 
\]
The above equality follows from an analogue of (\ref{3eq21}), i.e., 
\[
(\nabla _{v}^{\Gamma }F)(\gamma )=\left\langle (\nabla ^{\Gamma }F)(\gamma
),v\right\rangle _{T_{\gamma }\Gamma }, 
\]
where $v\in V_{0}(X)$ and the directional derivative $\nabla _{v}^{\Gamma }$
is as in Definition \ref{3eq9} with $\omega $ replaced by $\gamma $.
Moreover for the the adjoint operator to the gradient $\nabla ^{\Gamma }$ on 
$L^{2}(\Gamma ,\pi _{\sigma })$ equation (\ref{3eq16}) becomes equal to 
\begin{eqnarray*}
&&\int_{\Gamma }(\nabla _{v}^{\Gamma }F)(\gamma )G(\gamma )d\pi _{\sigma
}(\gamma ) \\
&=&-\int_{\Gamma }F(\gamma )(\nabla _{v}^{\Gamma }G)(\gamma )d\pi _{\sigma
}(\gamma )-\int_{\Gamma }F(\gamma )G(\gamma )B_{v}^{\pi _{\sigma }}(\gamma
)d\pi _{\sigma }(\gamma ),
\end{eqnarray*}
or 
\[
\nabla _{v}^{\Gamma *}=-\nabla _{v}^{\Gamma }-B_{v}^{\pi _{\sigma }} 
\]
as an operator equality on the domain $\mathcal{F}C_{b}^{\infty }(\mathcal{D}%
,\Gamma )$ in $L^{2}(\Gamma ,\pi _{\sigma })$ and $B_{v}^{\pi _{\sigma }}$
stands for the logarithmic derivative of the Poisson measure $\pi _{\sigma }$
along $v$, i.e., 
\[
B_{v}^{\pi _{\sigma }}(\gamma ):=\langle \gamma ,\beta _{v}^{\sigma }\rangle
=\int_{X}[\langle \beta ^{\sigma }(x),v(x)\rangle _{T_{x}X}+\mathrm{div}%
^{X}v(x)]d\gamma (x). 
\]

\subsection{A lifting of the geometry\label{3eq58}}

In the consideration above we have constructed some objects related to the
differential geometry of the space $\Omega $. Now we present an
interpretation of all the above formulas via a simple ``lifting rule''.

Any function $\varphi \in \mathcal{D}$ generates a (cylinder) function on $%
\Omega $ by the formula 
\begin{equation}
L_{\varphi }(\omega ):=\left\langle \omega ,\varphi \right\rangle ,\;\omega
\in \Omega .  \label{3eq27}
\end{equation}
We will call $L_{\varphi }$ the lifting of $\varphi $. As before any vector
field $v\in V_{0}(X)$ can be considered as a vector field on $\Omega $ (the
lifting of $v$) which we denote by $L_{v}$, see Definition \ref{3eq28}. For $%
v,w\in V_{0}(X)$ formula (\ref{3eq29}) can be written as 
\begin{equation}
\left\langle L_{v},L_{w}\right\rangle _{T_{\omega }\Omega }=L_{\left\langle
v,w\right\rangle _{TX}}(\omega ),  \label{3eq30}
\end{equation}
i.e., the scalar product of lifting vector fields is computed as the lifting
of the scalar product $\left\langle v(x),w(x)\right\rangle _{T_{x}X}=\varphi
(x)$. This rule can be used as a definition of the tangent space $T_{\omega
}\Omega $.

Formula (\ref{3eq10}) has now the following interpretation: 
\[
(\nabla _{v}^{\Omega }L_{\varphi })(\omega )=L_{\nabla _{v}^{X}\varphi
}(\omega ),\;\omega \in \Omega , 
\]
and the gradient of $L_{\varphi }$ is nothing but the lifting of the
corresponding underlying gradient: 
\[
(\nabla ^{\Omega }L_{\varphi })(\omega )=L_{\nabla ^{X}\varphi }(\omega ). 
\]

As follows from (\ref{3eq25}) the logarithmic derivative $B_{v}^{\pi
_{\sigma }^{\tau }}:\Omega \rightarrow \Bbb{R}$ is obtained via the same
lifting procedure of the corresponding logarithmic derivative $\beta
_{v}^{\sigma }:X\rightarrow \Bbb{R}$, namely, 
\[
B_{v}^{\pi _{\sigma }^{\tau }}(\omega )=L_{\beta _{v}^{\sigma }}(\gamma
_{\omega }). 
\]
Or, equivalently, one has for the divergence of a lifted vector field: 
\[
\mathrm{div}_{\pi _{\sigma }^{\tau }}^{\Omega }(L_{v})=L_{\mathrm{div}%
_{\sigma }^{X}v}. 
\]

\subsection{Representations of the Lie algebra of vector fields\label{3eq65}}

Using the property of quasi-invariance of the compound Poisson measure $\pi
_{\sigma }^{\tau }$ we can define a unitary representation of the
diffeomorphism group \textrm{Diff}$_{0}(X)$ in the space $L^{2}(\pi _{\sigma
}^{\tau })$, see \cite{GGV75}. Namely, for $\phi \in $\textrm{Diff}$_{0}(X)$
we define a unitary operator 
\[
(V_{\pi _{\sigma }^{\tau }}(\phi )F)(\omega ):=F(\phi (\omega ))\sqrt{\frac{%
d\pi _{\sigma }^{\tau }(\phi (\omega ))}{d\pi _{\sigma }^{\tau }(\omega )}}%
,\;F\in L^{2}(\pi _{\sigma }^{\tau }). 
\]
Then we have 
\[
V_{\pi _{\sigma }^{\tau }}(\phi _{1})V_{\pi _{\sigma }^{\tau }}(\phi
_{2})=V_{\pi _{\sigma }^{\tau }}(\phi _{1}\circ \phi _{2}),\;\phi _{1},\phi
_{2}\in \mathrm{Diff}_{0}(X). 
\]
as in Subsection \ref{3eq31}, to any vector field $v\in V_{0}(X)$ there
corresponds a one-parameter subgroup of diffeomorphisms $\phi _{t}^{v}$,$%
\,t\in \Bbb{R}$. It generates a one-parameter unitary group 
\begin{equation}
V_{\pi _{\sigma }^{\tau }}(\phi _{t}^{v}):=\exp [itJ_{\pi _{\sigma }^{\tau
}}(v)],\;t\in \Bbb{R,}  \label{3eq41}
\end{equation}
where $J_{\pi _{\sigma }^{\tau }}(v)$ denotes the self-adjoint generator of
this group.

\begin{proposition}
For any $v\in V_{0}(X)$ the following operator equality on the domain $%
\mathcal{F}C_{b}^{\infty }(\mathcal{D},\Omega )$ holds: 
\begin{equation}
J_{\pi _{\sigma }^{\tau }}(v)=\frac{1}{i}\nabla _{v}^{\Omega }+\frac{1}{2i}%
B_{v}^{\pi _{\sigma }^{\tau }}.  \label{3eq32}
\end{equation}
\end{proposition}

\noindent \textbf{Proof.} Let $F\in \mathcal{F}C_{b}^{\infty }(\mathcal{D}%
,\Omega )$ be given. Then differentiating the left hand side of (\ref{3eq41}%
) at $t=0$ we get 
\begin{eqnarray*}
\frac{d}{dt}(V_{\pi _{\sigma }^{\tau }}(\phi _{t}^{v})F)(\omega )|_{t=0} &=&%
\frac{d}{dt}F(\phi _{t}^{v}(\omega ))|_{t=0}+F(\omega )\frac{1}{2}\frac{d}{dt%
}\left. \frac{d\pi _{\sigma }^{\tau }(\phi _{t}^{v}(\omega ))}{d\pi _{\sigma
}^{\tau }(\omega )}\right| _{t=0} \\
&=&(\nabla _{v}^{\Omega }F)(\omega )+\frac{1}{2}F(\omega )B_{v}^{\pi
_{\sigma }^{\tau }}(\omega ),
\end{eqnarray*}
where we have used the form of the operator $V_{\pi _{\sigma }^{\tau }}(\phi
_{t}^{v})$, the definition of the directional derivative $\nabla
_{v}^{\Omega }$ and Theorem \ref{3eq42}. On the other hand the same
procedure on the right hand side of (\ref{3eq41}) produce $i(J_{\pi _{\sigma
}^{\tau }}(v)F)(\omega )$. Hence the result of the proposition
follows.\hfill $\blacksquare $

\begin{remark}
More generally, one can study a family of self-adjoint operators $J(v)$, $%
v\in V_{0}(X)$, in a Hilbert space $\mathcal{H}$ which gives a
representation of the Lie algebra $V_{0}(X)$ in the sense of the following
commutation relation: 
\begin{equation}
[J(v_{1}),J(v_{2})]=-iJ([v_{1},v_{2}])  \label{3eq33}
\end{equation}
(on a dense domain in $\mathcal{H}$), where $[v_{1},v_{2}]=\left\langle
v_{1},\nabla v_{2}\right\rangle _{TX}-\left\langle v_{2},\nabla
v_{1}\right\rangle _{TX}$ is the Lie-bracket of the vector fields $%
v_{1},v_{2}\in V_{0}(X)$. In the case discussed, this relation is a direct
consequence of (\ref{3eq32}). Thus, we have constructed a compound Poisson
space representation of the Lie algebra $V_{0}(X)$.
\end{remark}

Let us define, in addition, a unitary representation of the additive group $%
\mathcal{D}$ given by the formula 
\[
(U_{\pi _{\sigma }^{\tau }}(f)F)(\omega ):=\exp (i\left\langle \omega
,f\right\rangle )F(\omega ),\;F\in L^{2}(\pi _{\sigma }^{\tau }),\;\omega
\in \Omega , 
\]
for any $f\in \mathcal{D}$. As usual, the semi-direct product $\mathcal{G}:=%
\mathcal{D}\wedge \mathrm{Diff}_{0}(X)$ of the groups $\mathcal{D}$ and 
\textrm{Diff}$_{0}(X)$ is defined as the set of pairs $(f,\phi )$ with
multiplication operation 
\[
(f_{1},\phi _{1})(f_{2},\phi _{2})=(f_{1}+f_{2}\circ \phi _{1},\phi
_{2}\circ \phi _{1}), 
\]
see e.g. \cite{GGV75}. Let us introduce for any element $(f,\phi )\in 
\mathcal{G}$ the following operator on $L^{2}(\pi _{\sigma }^{\tau }):$%
\[
W_{\pi _{\sigma }^{\tau }}(f,\phi ):=U_{\pi _{\sigma }^{\tau }}(f)V_{\pi
_{\sigma }^{\tau }}(\phi ). 
\]
These operators are unitary and form a representation of the group $\mathcal{%
G}$. If we introduce multiplication operators $\rho _{\pi _{\sigma }^{\tau
}}(f)$, $f\in \mathcal{D}$, as self-adjoint operators on $L^{2}(\pi _{\sigma
}^{\tau })$ which are defined for $F\in \mathcal{F}C_{b}^{\infty }(\mathcal{D%
},\Omega )$ by the formula 
\[
(\rho _{\pi _{\sigma }^{\tau }}(f)F)(\omega ):=\left\langle \omega
,f\right\rangle F(\omega ),\;\omega \in \Omega , 
\]
then $U_{\pi _{\sigma }^{\tau }}(f)=\exp [i\rho _{\pi _{\sigma }^{\tau
}}(f)] $ and the form of the multiplication in $\mathcal{G}$ implies 
\[
\lbrack \rho _{\pi _{\sigma }^{\tau }}(f),J_{\pi _{\sigma }^{\tau
}}(v)]=i\rho _{\pi _{\sigma }^{\tau }}(\nabla _{v}^{X}f) 
\]
(on a dense domain in $L^{2}(\pi _{\sigma }^{\tau })$) for all $f\in 
\mathcal{D}$, $v\in \mathrm{Diff}_{0}(X)$. We also have the relation $[\rho
_{\pi _{\sigma }^{\tau }}(f_{1}),\rho _{\pi _{\sigma }^{\tau }}(f_{2})]=0$.
The family of operators $J_{\pi _{\sigma }^{\tau }}(v),\rho _{\pi _{\sigma
}^{\tau }}(f)$, $v\in V_{0}(X)$, $f\in \mathcal{D}$, thus forms a compound
Poisson representation of an infinite-dimensional Lie algebra. In the
particular case when $\tau =\varepsilon _{1}$ this representation is known
as Lie algebra of currents in non relativistic quantum field theory, e.g.~%
\cite{GGPS74}.

\section{Intrinsic Dirichlet forms on compound Poisson space\label{3eq60}}

\subsection{Definition of the intrinsic Dirichlet form\label{3eq62}}

We start with introducing some useful spaces of cylinder functions on $%
\Omega $ in addition to $\mathcal{F}C_{b}^{\infty }(\mathcal{D},\Omega )$.
By $\mathcal{F}\mathcal{P}(\mathcal{D},\Omega )$ we denote the set of all
cylinder functions of the form (\ref{3eq8}) in which the generating function $%
g_{F}$ is a polynomial on $\Bbb{R}^{N}$, i.e., $g_{F}\in \mathcal{P}(\Bbb{R}%
^{N})$. Analogously we define $\mathcal{F}C_{p}^{\infty }(\mathcal{D},\Omega
)$, where now $g_{F}\in C_{p}^{\infty }(\Bbb{R}^{N})$ (the set of all $%
C^{\infty }$-functions $f$ on $\Bbb{R}^{N}$ such that $f$ and all its
partial derivatives of any order are polynomially bounded).

We have obviously 
\begin{eqnarray*}
\mathcal{F}C_{b}^{\infty }(\mathcal{D},\Omega ) &\subset &\mathcal{F}%
C_{p}^{\infty }(\mathcal{D},\Omega ), \\
\mathcal{F}\mathcal{P}(\mathcal{D},\Omega ) &\subset &\mathcal{F}%
C_{p}^{\infty }(\mathcal{D},\Omega ),
\end{eqnarray*}
and these spaces are algebras with respect to the usual operations. The
existence of the Laplace transform $l_{\pi _{\sigma }^{\tau }}(f)$, $f\in 
\mathcal{D}$, implies $\mathcal{F}C_{p}^{\infty }(\mathcal{D},\Omega
)\subset L^{2}(\pi _{\sigma }^{\tau })$.

Note that after the embedding $\Omega \hookrightarrow \mathcal{D}^{\prime }$
(see Subsection \ref{3eq22}) and a natural extension to $\mathcal{D}^{\prime
} $ the space $\mathcal{F}\mathcal{P}(\mathcal{D},\mathcal{D}^{\prime })$ is
nothing but the well-known space of cylinder polynomials on $\mathcal{D}%
^{\prime }$, see \cite[Chap.~2]{BK88}.

\begin{definition}
For $F,G\in \mathcal{F}C_{p}^{\infty }(\mathcal{D},\Omega )$ we introduce a
pre-Dirichlet form as 
\begin{equation}
\mathcal{E}_{\pi _{\sigma }^{\tau }}^{\Omega }(F,G)=\int_{\Omega }\langle
(\nabla ^{\Omega }F)(\omega ),(\nabla ^{\Omega }G)(\omega )\rangle
_{T_{\omega }\Omega }d\pi _{\sigma }^{\tau }(\omega ).  \label{3eq34}
\end{equation}
\end{definition}

Note that for $F,G\in \mathcal{F}C_{p}^{\infty }(\mathcal{D},\Omega )$
formula (\ref{3eq11}) is still valid and therefore 
\[
\left\langle \nabla ^{\Omega }F,\nabla ^{\Omega }G\right\rangle _{T\Omega
}\in \mathcal{F}C_{p}^{\infty }(\mathcal{D},\Omega ), 
\]
such that (\ref{3eq34}) is well-defined.

We will call $\mathcal{E}_{\pi _{\sigma }^{\tau }}^{\Omega }$ the intrinsic
pre-Dirichlet form corresponding to the compound Poisson measure $\pi
_{\sigma }^{\tau }$ on $\Omega $. The name ``intrinsic'' means that $%
\mathcal{E}_{\pi _{\sigma }^{\tau }}^{\Omega }$ is associated with the
geometry of $\Omega $ generated by the original Riemannian structure of $X$,
in particular, by the intrinsic gradient $\nabla ^{\Omega }$. In the next
subsection we shall prove the closability of $\mathcal{E}_{\pi _{\sigma
}^{\tau }}^{\Omega }$.

\subsection{Intrinsic Dirichlet operators\label{3eq63}}

Let us introduce a differential operator $H_{\pi _{\sigma }^{\tau }}^{\Omega
}$ on the domain $\mathcal{F}C_{b}^{\infty }(\mathcal{D},\Omega )$ which is
given on any $F\in \mathcal{F}C_{b}^{\infty }(\mathcal{D},\Omega )$ of the
form 
\begin{equation}
F(\omega )=g_{F}(\langle \omega ,\varphi _{1}\rangle ,\ldots ,\langle \omega
,\varphi _{N}\rangle ),\;\omega \in \Omega ,g_{F}\in C_{b}^{\infty }(\Bbb{R}%
^{N}),\varphi _{1},\ldots ,\varphi _{N}\in \mathcal{D},  \label{3eq35}
\end{equation}
by the formula 
\begin{eqnarray}
&&(H_{\pi _{\sigma }^{\tau }}^{\Omega }F)(\omega )  \nonumber \\
&\mbox{$:=$}&\sum_{i,j=1}^{N}\frac{\partial ^{2}g_{F}}{\partial
s_{i}\partial s_{j}}(\langle \omega ,\varphi _{1}\rangle ,\ldots ,\langle
\omega ,\varphi _{N}\rangle )\int_{X}\langle \nabla ^{X}\varphi
_{i}(x),\nabla ^{X}\varphi _{j}(x)\rangle _{T_{x}X}d\omega (x)  \nonumber \\
&&-\sum_{i=1}^{N}\frac{\partial g_{F}}{\partial s_{i}}(\langle \omega
,\varphi _{1}\rangle ,\ldots ,\langle \omega ,\varphi _{N}\rangle
)\int_{X}\triangle ^{X}\varphi _{i}(x)d\omega (x)  \label{3eq38} \\
&&-\sum_{i=1}^{N}\frac{\partial g_{F}}{\partial s_{i}}(\langle \omega
,\varphi _{1}\rangle ,\ldots ,\langle \omega ,\varphi _{N}\rangle
)\int_{X}\langle \nabla ^{X}\varphi _{i}(x),\beta ^{\sigma }(x)\rangle
_{T_{x}X}d\omega (x),  \nonumber
\end{eqnarray}
where $\triangle ^{X}$ denotes the Laplace-Beltrami operator on $X$. In this
formula all expressions are from $\mathcal{F}C_{b}^{\infty }(\mathcal{D}%
,\Omega )$ or have the form $\left\langle \omega ,\psi \right\rangle $, $%
\omega \in \Omega $, $\psi \in \mathcal{D}$, except for the functions $%
\left\langle \omega ,h_{i}\right\rangle $, $\omega \in \Omega $, with 
\[
h_{i}(x)=\left\langle (\nabla ^{X}\varphi _{i})(x),\beta ^{\sigma
}(x)\right\rangle _{T_{x}X},\;x\in X,i=1,\ldots ,N. 
\]

To clarify the situation with these functions note that due to the
assumption on $\sigma $ we have $\rho ^{1/2}\in H_{loc}^{1,2}(X)$ which
gives $h_{i}\in L^{1}(\sigma )$ and these functions have compact supports.
Therefore $h_{i}\in L^{1}(\sigma )$, $j=1,\ldots ,N$. On the other hand we
know that a function $\left\langle \omega ,f\right\rangle $, $\omega \in
\Omega $, is from $L^{2}(\pi _{\sigma }^{\tau })$ if $f\in L^{1}(\sigma
)\cap L^{2}(\sigma )$. The latter follows from the formula for the second
moment of the measure $\pi _{\sigma }^{\tau }$, namely 
\begin{equation}
\int_{\Omega }\left\langle \omega ,f\right\rangle ^{2}d\pi _{\sigma }^{\tau
}(\omega )=m_{2}(\tau )\int_{X}f^{2}(x)d\sigma (x)+(m_{1}(\tau ))^{2}\left(
\int_{X}f(x)d\sigma (x)\right) ^{2},  \label{3eq36}
\end{equation}
where $m_{1}(\tau )$ and $m_{2}(\tau )$ are the first and second moment of
the measure $\tau $ on $\Bbb{R}^{+}$, respectively. Equation (\ref{3eq36}) is
a direct consequence of (\ref{3eq37}). As a result the right hand side of (%
\ref{3eq38}) is well-defined. To show that the operator $H_{\pi _{\sigma
}^{\tau }}^{\Omega }$ is well-defined we still have to show that its
definition does not depend on the representation of $F$ in (\ref{3eq35})
which will be done below.

\begin{remark}
In the applications to the study of unitary representations of the group 
\emph{Diff}$_{0}(X)$ given by compound Poisson measures, there is usually an
additional assumption on the smoothness of the density $\rho :=d\sigma /dm$,
namely $\rho \in C^{\infty }(X)$, $\rho (x)>0$, $x\in X$, see e.g.~\cite
{GGV75}. In this case it is obvious that the operator $H_{\pi _{\sigma
}^{\tau }}^{\Omega }$ preserves the spaces $\mathcal{F}C_{p}^{\infty }(%
\mathcal{D},\Omega )$ and $\mathcal{F}\mathcal{P}(\mathcal{D},\Omega )$.
\end{remark}

Let us also consider the classical pre-Dirichlet form corresponding to the
measure $\sigma $ on $X:$%
\begin{equation}
\mathcal{E}_{\sigma }^{X}(\varphi ,\psi ):=\int_{X}\langle \nabla
^{X}\varphi (x),\nabla ^{X}\psi (x)\rangle _{T_{x}X}d\sigma (x),
\label{3eq39}
\end{equation}
where $\varphi ,\psi \in \mathcal{D}$. This form is associated with the
Dirichlet operator $H_{\sigma }^{X}$ which is given on $\mathcal{D}$ by 
\begin{equation}
(H_{\sigma }^{X}\varphi )(x):=-\triangle ^{X}\varphi (x)-\left\langle \beta
^{\sigma }(x),\nabla ^{X}\varphi (x)\right\rangle _{T_{x}X},  \label{3eq40}
\end{equation}
and which satisfies 
\[
\mathcal{E}_{\sigma }^{X}(\varphi ,\psi )=(H_{\sigma }^{X}\varphi ,\psi
)_{L^{2}(\sigma )},\;\varphi ,\psi \in \mathcal{D}. 
\]

The closure of this form on $L^{2}(\sigma )$ is defined by $(\mathcal{E}%
_{\sigma }^{X},D(\mathcal{E}_{\sigma }^{X}))$. Note that $D(\mathcal{E}%
_{\sigma }^{X})$ is nothing but the Sobolev space of order 1 in $%
L^{2}(\sigma )$ (sometimes also denoted by $H_{0}^{1,2}(X,\sigma )$). $(%
\mathcal{E}_{\sigma }^{X},D(\mathcal{E}_{\sigma }^{X}))$ generates a
positive self-adjoint operator in $L^{2}(\sigma )$ (the so-called
Friedrich's extension of $H_{\sigma }^{X}$, see e.g.~\cite{BKR97} and \cite
{RS75}).

For this extension we preserve the previous notation $H_{\sigma }^{X}$ and
denote the domain by $D(H_{\sigma }^{X})$. Using the underlying Dirichlet
operator we obtain the representation 
\begin{eqnarray*}
&&(H_{\pi _{\sigma }^{\tau }}^{\Omega }F)(\omega ) \\
&=&\sum_{i,j=1}^{N}\frac{\partial ^{2}g_{F}}{\partial s_{i}\partial s_{j}}%
(\langle \omega ,\varphi _{1}\rangle ,\ldots ,\langle \omega ,\varphi
_{N}\rangle )\langle \nabla ^{X}\varphi _{i},\nabla ^{X}\varphi _{j}\rangle
_{T_{\omega }\Omega } \\
&&+\sum_{i=1}^{N}\frac{\partial g_{F}}{\partial s_{i}}(\langle \omega
,\varphi _{1}\rangle ,\ldots ,\langle \omega ,\varphi _{N}\rangle )\langle
\omega ,H_{\sigma }^{X}\varphi _{i}\rangle .
\end{eqnarray*}

Let us define for any $\omega \in \Omega $, $x,y\in X$%
\begin{eqnarray*}
&&(\nabla ^{\Omega }\nabla ^{\Omega }F)(\omega ,x,y) \\
&\mbox{$:=$}&\sum_{i,j=1}^{N}\frac{\partial ^{2}g_{F}}{\partial
s_{i}\partial s_{j}}(\langle \omega ,\varphi _{1}\rangle ,\ldots ,\langle
\omega ,\varphi _{N}\rangle )\nabla ^{X}\varphi _{i}(x)\otimes \nabla
^{X}\varphi _{j}(y)\in T_{\omega }\Omega \otimes T_{\omega }\Omega .
\end{eqnarray*}
Then 
\begin{eqnarray*}
\triangle ^{\Omega }F(\omega ) &\mbox{$:=$}&\mathrm{Tr}(\nabla ^{\Omega
}\nabla ^{\Omega }F)(\omega ) \\
&=&\sum_{i,j=1}^{N}\frac{\partial ^{2}g_{F}}{\partial s_{i}\partial s_{j}}%
(\langle \omega ,\varphi _{1}\rangle ,\ldots ,\langle \omega ,\varphi
_{N}\rangle )\langle \nabla ^{X}\varphi _{i},\nabla ^{X}\varphi _{j}\rangle
_{T_{\omega }\Omega }.
\end{eqnarray*}
Hence the operator $H_{\pi _{\sigma }^{\tau }}^{\Omega }$ can be written as 
\[
(H_{\pi _{\sigma }^{\tau }}^{\Omega }F)(\omega )=-(\triangle ^{\Omega
}F)(\omega )-\langle \omega ,\mathrm{div}_{\sigma }^{X}(\nabla ^{\Omega
}F)(\omega ;\cdot )\rangle .
\]

The following theorem implies that both $H_{\pi _{\sigma }^{\tau }}^{\Omega
} $ and $\triangle ^{\Omega }$ are well-defined as linear operators on $%
\mathcal{F}C_{b}^{\infty }(\mathcal{D},\Omega )$, i.e., independently of the
representation of $F$ on (\ref{3eq35}).

\begin{theorem}
\label{3eq45}The operator $H_{\pi _{\sigma }^{\tau }}^{\Omega }$ is
associated with the intrinsic Dirichlet form $\mathcal{E}_{\pi _{\sigma
}^{\tau }}^{\Omega }$, i.e., for all $F,G\in \mathcal{F}C_{b}^{\infty }(%
\mathcal{D},\Omega )$ we have 
\[
\mathcal{E}_{\pi _{\sigma }^{\tau }}^{\Omega }(F,G)=(H_{\pi _{\sigma }^{\tau
}}^{\Omega }F,G)_{L^{2}(\pi _{\sigma }^{\tau })}, 
\]
or 
\[
H_{\pi _{\sigma }^{\tau }}^{\Omega }=-\mathrm{div}_{\pi _{\sigma }^{\tau
}}^{\Omega }\nabla ^{\Omega }\;\mathrm{on\;}\mathcal{F}C_{b}^{\infty }(%
\mathcal{D},\Omega ). 
\]
We call $H_{\pi _{\sigma }^{\tau }}^{\Omega }$ the intrinsic Dirichlet
operator of the measure $\pi _{\sigma }^{\tau }$.
\end{theorem}

\noindent \textbf{Proof.} For any $F\in \mathcal{F}C_{b}^{\infty }(\mathcal{D%
},\Omega )$ of the form (\ref{3eq35}) we have 
\[
(\nabla ^{\Omega }F)(\omega ;x)=\sum_{j=1}^{N}\frac{\partial g_{F}}{\partial
s_{i}}(\left\langle \omega ,\varphi _{1}\right\rangle ,\ldots ,\left\langle
\omega ,\varphi _{N}\right\rangle )\nabla ^{X}\varphi _{i}(x). 
\]
By (\ref{3eq43}) we conclude that 
\[
\mathrm{div}_{\pi _{\sigma }^{\tau }}^{\Omega }(\nabla ^{\Omega }F)=-H_{\pi
_{\sigma }^{\tau }}^{\Omega }F 
\]
which by (\ref{3eq44}) for $F,G\in \mathcal{F}C_{b}^{\infty }(\mathcal{D}%
,\Omega )$ gives 
\begin{eqnarray*}
(H_{\pi _{\sigma }^{\tau }}^{\Omega }F,G)_{L^{2}(\pi _{\sigma }^{\tau })}
&=&-\int_{\Omega }\mathrm{div}_{\pi _{\sigma }^{\tau }}^{\Omega }(\nabla
^{\Omega }F)(\omega )G(\omega )d\pi _{\sigma }^{\tau }(\omega ) \\
&=&\int_{\Omega }\left\langle (\nabla ^{\Omega }F)(\omega ),(\nabla ^{\Omega
}G)(\omega )\right\rangle _{T_{\omega }\Omega }d\pi _{\sigma }^{\tau
}(\omega ).
\end{eqnarray*}
\hfill $\blacksquare $

\begin{remark}
\begin{enumerate}
\item  In the case $\sigma =m$ and $\tau =\varepsilon _{1}$ we call the
Dirichlet form $\mathcal{E}_{\pi _{\sigma }}^{\Gamma }$ the canonical
Dirichlet form on $\Gamma $. The canonical Dirichlet form and canonical
Dirichlet operator $H_{\pi _{\sigma }}^{\Gamma }$ are defined directly in
terms of the Riemannian geometry of $X$.

\item  The operator $H_{\pi _{\sigma }}^{\Gamma }$ can be naturally extended
to cylinder functions of the form 
\[
F(\omega ):=\exp (\langle \omega ,\varphi \rangle ),\varphi \in \mathcal{D}%
,\omega \in \Omega , 
\]
since such $F$ belongs to $L^{2}(\Gamma ,\pi _{\sigma })$. We then have 
\[
(H_{\pi _{\sigma }}^{\Gamma }\exp (\langle \cdot ,\varphi \rangle ))(\gamma
)=\langle \omega ,H_{\sigma }^{X}\varphi -\left| \nabla ^{X}\varphi \right|
_{TX}^{2}\rangle \exp (\langle \gamma ,\varphi \rangle ) 
\]
and for $\sigma =m$ and $\tau =\varepsilon _{1}$%
\[
(H_{\pi _{m}}^{\Gamma }\exp (\langle \cdot ,\varphi \rangle ))(\gamma
)=-\langle \triangle ^{X}\varphi +\left| \nabla ^{X}\varphi \right|
_{TX}^{2}\rangle \exp (\langle \gamma ,\varphi \rangle ). 
\]
\end{enumerate}
\end{remark}

As an immediate consequence of Theorem \ref{3eq45} we obtain

\begin{corollary}
$(\mathcal{E}_{\pi _{\sigma }^{\tau }}^{\Omega },\mathcal{F}C_{b}^{\infty }(%
\mathcal{D},\Omega ))$ is closable on $L^{2}(\pi _{\sigma }^{\tau })$. Its
closure $(\mathcal{E}_{\pi _{\sigma }^{\tau }}^{\Omega },$\linebreak $D(%
\mathcal{E}_{\pi _{\sigma }^{\tau }}^{\Omega }))$ is associated with a
positive definite self-adjoint operator, the Fried\-richs extension of $%
H_{\pi _{\sigma }^{\tau }}^{\Omega }$ which we also denote by $H_{\pi
_{\sigma }^{\tau }}^{\Omega }$ (and its domain by $D(H_{\pi _{\sigma }^{\tau
}}^{\Omega })$).
\end{corollary}

Clearly, $\nabla ^{\Omega }$ also extends to $D(\mathcal{E}_{\pi _{\sigma
}^{\tau }}^{\Omega })$. We denote this extension again by $\nabla ^{\Omega }$%
.

\begin{corollary}
Let 
\begin{eqnarray*}
F(\omega )=g_{F}(\langle \omega ,\varphi _{1}\rangle ,\ldots ,\langle \omega
,\varphi _{N}\rangle ), &&\omega \in \Omega ,g_{F}\in C_{b}^{\infty }(\Bbb{R}%
^{N}), \\
&&\varphi _{1},\ldots ,\varphi _{N}\in D(\mathcal{E}_{\sigma }^{\Omega }).
\end{eqnarray*}
Then $F\in D(\mathcal{E}_{\pi _{\sigma }^{\tau }}^{\Omega })$ and 
\[
(\nabla ^{\Omega }F)(\cdot )=\sum_{i=1}^{N}\frac{\partial g_{F}}{\partial
s_{i}}(\langle \cdot ,\varphi _{1}\rangle ,\ldots ,\langle \cdot ,\varphi
_{N}\rangle )\nabla ^{X}\varphi _{i}. 
\]
\end{corollary}

\noindent \textbf{Proof.} By approximation this is an immediate consequence
of (\ref{3eq11}) and the fact that for all $1\leq i\leq N$%
\[
(m_{1}(\tau ))^{-1}\int_{\Omega }\langle \omega ,|\nabla ^{X}\varphi
_{i}|_{TX}^{2}\rangle d\pi _{\sigma }^{\tau }(\omega )=\mathcal{E}_{\sigma
}^{X}(\varphi _{i},\varphi _{i}). 
\]
\hfill $\blacksquare $

\begin{remark}
Of course the domain of $\mathcal{E}_{\pi _{\sigma }^{\tau }}^{\Omega }$ $D(%
\mathcal{E}_{\pi _{\sigma }^{\tau }}^{\Omega })$ is nothing but the Sobolev
space $H_{0}^{1,2}(\Omega ,\pi _{\sigma }^{\tau })$ on $\Omega $ of order 1
in $L^{2}(\Omega ,\pi _{\sigma }^{\tau })$.
\end{remark}

\section{Identification of the process on compound Poisson space\label{3eq64}}

In this section we will prove the existence of a diffusion process
corresponding to our Dirichlet form $(\mathcal{E}_{\pi _{\sigma }^{\tau
}}^{\Omega },\mathcal{D}(\mathcal{E}_{\pi _{\sigma }^{\tau }}^{\Omega }))$.
For a general theory of processes corresponding to Dirichlet forms we refer
to \cite[Chap.~IV]{MR92}, see also \cite{F80}.

After all our preparation and taking into account the general description of
compound Poisson space given in Section \ref{3eq56} we will see that this
process is nothing but a direct consequence (``lifting'') of the
corresponding process on the space of simple configurations $\Gamma _{X}$,
see \cite[Sect.~6]{AKR97} for a detailed description. Let us clarify this in
more detail. After Section \ref{3eq56} there is an obvious identification
between compound configurations $\omega \in \Omega $ and a marked
configurations $(\gamma _{\omega },m_{\omega })\in \Omega _{X}^{\Bbb{R}_{+}}$
which gives the possibility to obtain an embedding from $L^{2}(\Gamma
_{X},\pi _{\sigma })$ into $L^{2}(\Omega _{X},\pi _{\sigma }^{\tau })$,
i.e., 
\[
L^{2}(\Omega _{X},\pi _{\sigma }^{\tau })\ni F(\omega )=G(\gamma _{\omega
})\in L^{2}(\Gamma _{X},\pi _{\sigma }).
\]
Hence all operators acting in $L^{2}(\Gamma _{X},\pi _{\sigma })$, e.g., $%
\nabla ^{\Gamma }$, $\nabla ^{\Gamma *}$, $H_{\pi _{\sigma }}^{\Gamma }$
etc.~are applicable on $L^{2}(\Omega _{X},\pi _{\sigma }^{\tau })$
w.r.t.~part of the variables. Moreover we have the following relation 
\[
(\nabla ^{\Omega }F)(\omega )=(\nabla ^{\Gamma }F)((\gamma _{\omega
},m_{\omega })),\;\omega =(\gamma _{\omega },m_{\omega })\in \Omega ,
\]
for the intrinsic gradient, from which everything else follows, see below.

Let us consider a probability measure $\tau $ on $\Bbb{R}_{+}$ (or more
general a probability measure on the space of marks $M$, cf.~Section \ref
{3eq56}). In what follows we always identify any compound configuration $%
\omega \in \Omega $ (or in general marked configuration) with $(\gamma
_{\omega },m_{\omega })$, i.e., 
\[
\Omega \ni \omega \rightsquigarrow (\gamma _{\omega },m_{\omega })
\]
and by Subsection \ref{3eq71} we have for any diffeomorphism $\phi \in 
\mathrm{Diff}_{0}(X)$ its action on $(\gamma _{\omega },m_{\omega })$ is
given by 
\[
\phi (\gamma _{\omega },m_{\omega })=(\phi (\gamma _{\omega }),m_{\omega }).
\]

It follows from Proposition \ref{3eq17} and the assumption on $\tau $ that
the Radon-Nikodym density of $\pi _{\sigma }^{\tau }$ and $\pi _{\sigma }$
with respect to the group \textrm{Diff}$_{0}(X)$ are equal, i.e., $p_{\phi
}^{\pi _{\sigma }^{\tau }}(\omega )=p_{\phi }^{\pi _{\sigma }}(\gamma
_{\omega })$, where $\gamma _{\omega }$ corresponds to $\omega $.

Let us compute the action of the gradient $\nabla ^{\Gamma }$ on cylinder
functions $F\in \mathcal{F}C_{b}^{\infty }(\mathcal{D},\Omega )$. To this
end let $v\in V_{0}(X)$ be a vector field on $X$ with compact support and $%
\phi _{t}^{v}$ the corresponding flow. Then by definition we have 
\[
(\nabla _{v}^{\Gamma }F)((\gamma _{\omega },m_{\omega })):=\frac{d}{dt}%
F((\phi _{t}^{v}(\gamma _{\omega }),m_{\omega }))|_{t=0}. 
\]
On cylinder functions of the form 
\[
F((\gamma _{\omega },m_{\omega }))=g_{F}(\langle (\gamma _{\omega
},m_{\omega }),\varphi _{1}\rangle ,\ldots ,\langle (\gamma _{\omega
},m_{\omega }),\varphi _{N}\rangle ), 
\]
where $\omega =(\gamma _{\omega },m_{\omega })\in \Omega $, $g_{F}\in
C_{b}^{\infty }(\Bbb{R}^{N})$, and $\varphi _{1},\ldots \varphi _{N}\in 
\mathcal{D}$, the above definition gives 
\[
\sum_{i=1}^{N}\frac{\partial g_{F}}{\partial s_{i}}((\gamma _{\omega
},m_{\omega }),\varphi _{1}\rangle ,\ldots ,\langle (\gamma _{\omega
},m_{\omega }),\varphi _{N}\rangle )\langle (\gamma _{\omega },m_{\omega
}),\nabla _{v}^{X}\varphi _{i}\rangle . 
\]
Therefore the following equality on $\mathcal{F}C_{b}^{\infty }(\mathcal{D}%
,\Omega )$ (dense in $L^{2}(\pi _{\sigma }^{\tau })$) of the directional
derivatives holds 
\[
(\nabla _{v}^{\Omega }F)(\omega )=(\nabla _{v}^{\Gamma }F)((\gamma _{\omega
},m_{\omega })) 
\]
which implies the equality between the intrinsic gradients, i.e., 
\[
(\nabla ^{\Omega }F)(\omega )=(\nabla ^{\Gamma }F)((\gamma _{\omega
},m_{\omega })). 
\]

From these considerations on the intrinsic gradient we get a relation
between the Dirichlet forms, namely 
\begin{eqnarray*}
\mathcal{E}_{\pi _{\sigma }^{\tau }}^{\Omega }(F,G) &=&\int_{\Omega }\langle
(\nabla ^{\Omega }F)(\omega ),(\nabla ^{\Omega }G)(\omega )\rangle
_{T_{\omega }\Omega }d\pi _{\sigma }^{\tau }(\omega ) \\
&=&\int_{\Omega }\langle (\nabla ^{\Gamma }F)((\gamma _{\omega },m_{\omega
})),(\nabla ^{\Gamma }G)((\gamma _{\omega },m_{\omega }))\rangle _{T_{\gamma
_{\omega }}\Gamma }d\pi _{\sigma }^{\tau }(\omega ) \\
&=&\mathcal{E}_{\pi _{\sigma }^{\tau }}^{\Gamma }(F,G).
\end{eqnarray*}

On the other hand the above relation between Dirichlet forms allowed us to
derive easily the following relation for the intrinsic Dirichlet operators 
\[
(H_{\pi _{\sigma }^{\tau }}^{\Omega }F)(\omega )=(H_{\pi _{\sigma },\gamma
_{\omega }}^{\Gamma }F)(\gamma _{\omega },m_{\omega }),\;F\in \mathcal{F}%
C_{b}^{\infty }(\mathcal{D},\Omega ), 
\]
where $H_{\pi _{\sigma },\gamma _{\omega }}^{\Gamma }$ acts w.r.t.~the
variable $\gamma _{\omega }$. Of course the corresponding semigroups (whose
generators are $H_{\pi _{\sigma }^{\tau }}^{\Omega }$ and $H_{\pi _{\sigma
},\gamma _{\omega }}^{\Gamma }$) are related by 
\[
e^{-tH_{\pi _{\sigma }^{\tau }}^{\Omega }}=e^{-tH_{\pi _{\sigma },\gamma
_{\omega }}^{\Gamma }}\otimes \mathbf{1}_{m_{\cdot }},\;t>0. 
\]

From this it follows that the process, $\Xi _{t}$, $t\geq 0$, on compound
Poisson space (or in marked Poisson space) is nothing but the equilibrium
process $X_{t}^{\gamma _{\omega }}$ (distorted Brownian motion on $\Omega $)
together with marks of the corresponding configuration, i.e., 
\[
\Xi _{t}=\{X_{t}^{\gamma _{\omega }},m_{\omega }\}. 
\]
For more detailed description of properties of the process $X_{t}^{\gamma
_{\omega }}$ we refer to \cite[Sect.~6]{AKR97}.

\section{Marked Poisson measures\label{3eq56}}

In this section we present some general results on marked Poisson measures
which generalizes the compound Poisson measures introduced in Subsection \ref
{3eq22}. For more detailed information on marked Poisson processes we refer
to \cite[Chap.~6]{BL95}, \cite[Chap.~5]{K93}, \cite{MM91}, and references
therein.

\subsection{Definition and measurable structure}

In this subsection we define and describe the space of marked configurations
as well as its associated measurable structure. Before we recall the
definition of the configuration space over the Cartesian space $X\times M$
between a Riemannian manifold $X$ and a complete separable metric space $M$.

The configuration space $\Gamma _{X\times M}$ over the Cartesian product $%
X\times M$ is defined as the set of all locally finite subsets
(configurations) in $X\times M$: 
\[
\Gamma _{X\times M}:=\{\hat{\gamma }\subset X\times M||\hat{\gamma }\cap
K|<\infty \;\mathrm{for\;any\;compact\;}K\subset X\times M\}. 
\]
For any $\Lambda \subset X$ we sometimes use the shorthand $\gamma _{\Lambda
}$ for $\gamma \cap \Lambda $, for any $\gamma \in \Gamma _{X}$ and define 
\[
\Gamma _{\Lambda }:=\{\gamma \in \Gamma _{X}|\gamma \cap (X\backslash
\Lambda )=\emptyset \}. 
\]
For any $n\in \Bbb{N}_{0}$ and $\Lambda \in \mathcal{O}_{c}(X)$ we introduce
the space of $n$ points configurations as 
\[
\Gamma _{\Lambda }^{(n)}:=\{\gamma \in \Gamma _{\Lambda }||\gamma
|=n\},\;\Gamma _{\Lambda }^{(0)}:=\{\emptyset \}. 
\]

Let us now introduce the space of marked configurations which will plays the
same role as $\Gamma _{X}$ played for Poisson measure but now for the marked
Poisson measure, see Subsection \ref{3eq50} below. It is defined as 
\[
\Omega _{X}^{M}:=\{\omega =(\gamma _{\omega },m_{\omega })|\gamma _{\omega
}\in \Gamma _{X},m_{\omega }\in M^{\omega }\}. 
\]
Here $M^{\omega }$ stands for the set of all maps $\gamma _{\omega }\ni
x\mapsto m_{x}\in M$. We may also write the marked configuration space $%
\Omega _{X}^{M}$ as a subspace of $\Gamma _{X\times M}$ as follows 
\[
\Omega _{X}^{M}:=\{\omega =\{(x,m_{x})\}\subset \Gamma _{X\times
M}|\{x\}=\gamma _{\omega }\in \Gamma _{X},m_{x}\in M\}. 
\]
For any $\Lambda \subset X$ we define in an analogous way the set $\Omega
_{\Lambda }^{M}$, i.e., 
\[
\Omega _{\Lambda }^{M}:=\{\omega =(\gamma _{\omega },m_{\omega })|\gamma
_{\omega }\in \Gamma _{\Lambda },m_{\omega }\in M^{\omega }\} 
\]
and 
\[
\Omega _{\Lambda }^{M}:=\{\omega =\{(x,m_{x})\}\subset \Gamma _{\Lambda
\times M}|\{x\}=\gamma _{\omega }\in \Gamma _{\Lambda },m_{x}\in M\}. 
\]

In order to describe the $\sigma $-algebra $\mathcal{B}(\Omega _{X}^{M})$ we
proceed as follows. Let $\Lambda \in \mathcal{O}_{c}(X)$ and $n\in \Bbb{N}%
_{0}:=\Bbb{N}\cup \{0\}$ be given. We define an equivalent relation $\sim $
on $(\Lambda \times M)^{n}$ setting 
\[
((x_{1},m_{x_{1}}),(x_{2},m_{x_{2}}),\ldots ,(x_{n},m_{x_{n}}))\sim
((y_{1},m_{y_{1}}),(y_{2},m_{y_{2}}),\ldots ,(y_{n},m_{y_{n}})) 
\]
iff there exists a permutation $\pi \in \frak{S}_{n}$ (the group of
permutations of $n$ elements) such that 
\[
(x_{i},m_{x_{i}})=(y_{\pi (i)},m_{y_{\pi (i)}}),\;\forall i=1,\ldots ,n. 
\]
Hence we obtain the quotient space $(\Lambda \times M)^{n}/\frak{S}_{n}$ by
means of $\sim $. Then we introduce the subset of $(\Lambda \times M)^{n}/%
\frak{S}_{n}$, $\Omega _{\Lambda }^{M}(n)$, defined as follows 
\[
\Omega _{\Lambda }^{M}(n):=\{((x_{1},m_{x_{1}}),\ldots
,(x_{1},m_{x_{1}}))|x_{i}\in \Lambda ,x_{i}\neq x_{j},\,i\neq j,m_{x_{i}}\in
M\}, 
\]
or equivalently 
\[
\Omega _{\Lambda }^{M}(n):=\{(\gamma _{\omega },m_{\omega })|\gamma _{\omega
}\in \Gamma _{\Lambda }^{(n)},m_{\omega }\in M^{\omega }\},\;\Omega
_{\Lambda }^{M}(0):=\{\emptyset \} 
\]

The space $\Omega _{\Lambda }^{M}(n)$ is endowed with the relative metric
from $(\Lambda \times M)^{n}/\frak{S}_{n}$, i.e., 
\[
\delta ([x],[y])=\inf_{x^{\prime }\in [x],y^{\prime }\in [y]}d^{n}(x^{\prime
},y^{\prime }), 
\]
where $d^{n}$ is the metric defined on $(\Lambda \times M)^{n}$ driven from
the original metrics on $X$ and $M$. Therefore $\Omega _{\Lambda }^{M}(n)$
becomes a metrizable topological space.

It is obvious that 
\[
\Omega _{\Lambda }^{M}=\bigsqcup_{n=0}^{\infty }\Omega _{\Lambda }^{M}(n). 
\]
This space can be equipped with the topology of disjoint union of
topological spaces, namely, the strongest topology on $\Omega _{\Lambda
}^{M} $ such that all the embeddings 
\[
i_{n}:\Omega _{\Lambda }^{M}(n)\to \Omega _{\Lambda }^{M},\;n\in \Bbb{N}_{0} 
\]
are continuous. $\mathcal{B}(\Omega _{\Lambda }^{M})$ stands for the
corresponding Borel $\sigma $-algebra.

For any $\Lambda \in \mathcal{O}_{c}(X)$ there are natural restriction maps 
\[
p_{\Lambda }:\Omega _{X}^{M}\to \Omega _{\Lambda }^{M} 
\]
defined by 
\begin{equation}
p_{\Lambda }(\gamma _{\omega },m_{\omega })=(\gamma _{\omega }\cap \Lambda
,m_{\omega |\gamma _{\omega }\cap \Lambda })\in \Omega _{\Lambda
}^{M},\;(\gamma _{\omega },m_{\omega })\in \Omega _{X}^{M}.  \label{3eq51}
\end{equation}
The topology on $\Omega _{X}^{M}$ is defined as the weakest topology making
all the mappings $p_{\Lambda }$ continuous. The associated Borel $\sigma $%
-algebra is denoted by $\mathcal{B}(\Omega _{X}^{M})$.

\subsection{The projective limit}

Finally we want to show that $\Omega _{X}^{M}$ coincides with the projective
limit of the family of topological spaces $\{\Omega _{\Lambda }^{M}|\Lambda
\in \mathcal{O}_{c}(X)\}$. First we recall the definition of projective
limit of topological spaces, see e.g. \cite[Chap.~3]{BD68} and \cite[Chap.~2]
{S71}.

\begin{definition}
\label{3eq48}Let $\Lambda _{1},\Lambda _{2}\in \mathcal{O}_{c}(X)$ be given
with $\Lambda _{1}\subset \Lambda _{2}$. There are natural maps 
\[
p_{\Lambda _{2},\Lambda _{1}}:\Omega _{\Lambda _{2}}^{M}\to \Omega _{\Lambda
_{1}}^{M} 
\]
defined by 
\[
p_{\Lambda _{2},\Lambda _{1}}(\gamma _{\omega },m_{\omega })=(\gamma
_{\omega }\cap _{\Lambda _{1}},m_{\omega |\gamma _{\omega }\cap \Lambda
_{1}})\in \Omega _{\Lambda _{1}}^{M},\;(\gamma _{\omega },m_{\omega })\in
\Omega _{\Lambda _{2}}^{M}. 
\]
The projective limit of the family $\{\Omega _{\Lambda }^{M}|\Lambda \in 
\mathcal{O}_{c}(X)\}$ denoted by 
\[
\limfunc{prlim}_{\Lambda \in \mathcal{O}_{c}(X)}\Omega _{\Lambda }^{M} 
\]
is a topological space $\Omega $ and a family of continuous projections 
\[
P_{\Lambda }:\Omega \to \Omega _{\Lambda }^{M},\;\Lambda \in \mathcal{O}%
_{c}(X), 
\]
such that the following two conditions are satisfied:

\begin{enumerate}
\item  \label{3eq47}If $\Lambda _{1},\Lambda _{2}\in \mathcal{O}_{c}(X)$
with $\Lambda _{1}\subset \Lambda _{2}$, then 
\[
P_{\Lambda _{1}}=p_{\Lambda _{2},\Lambda _{1}}\circ P_{\Lambda _{2}}. 
\]

\item  \label{3eq49}If $\Omega ^{\prime }$ is a topological space and 
\[
P_{\Lambda }^{\prime }:\Omega ^{\prime }\to \Omega _{\Lambda }^{M},\;\Lambda
\in \mathcal{O}_{c}(X), 
\]
a family of continuous projections which fulfills Condition \ref{3eq47}
above, then there exists a unique continuous map $u:\Omega ^{\prime
}\rightarrow \Omega $ such that $P_{\Lambda }^{\prime }=P_{\Lambda }\circ u$%
, for all $\Lambda \in \mathcal{O}_{c}(X)$.
\end{enumerate}
\end{definition}

\begin{remark}
The projective limit of the family $\{\Omega _{\Lambda }^{M}|\Lambda \in 
\mathcal{O}_{c}(X)\}$ exists and is unique in the following sense: let $%
\Omega $ and $\Omega ^{\prime }$ be projective limits, then there exists a
map $u:\Omega \to \Omega ^{\prime }$ such that $u$ and $u^{-1}$ are
continuous, see e.g. \cite[Chap.~4]{P67}.
\end{remark}

\begin{theorem}
The space of marked configurations $\Omega _{X}^{M}$ is the projective limit
of the family $\{\Omega _{\Lambda }^{M}|\Lambda \in \mathcal{O}_{c}(X)\}$
together with the family of projections $\{p_{\Lambda }|\Lambda \in \mathcal{%
O}_{c}(X)\}$ (defined in (\ref{3eq51})) and 
\[
\mathcal{B}(\Omega _{X}^{M})=\sigma (p_{\Lambda }^{-1}(\mathcal{B}(\Omega
_{\Lambda }^{M}));\Lambda \in \mathcal{O}_{c}(X)). 
\]
In other words there exists a bicontinuous bijective mapping between $\Omega
_{X}^{M}$ and the projective limit. This will be denoted by $\Omega
_{X}^{M}\simeq \limfunc{prlim}_{\Lambda \in \mathcal{O}_{c}(X)}\Omega
_{\Lambda }^{M}$.
\end{theorem}

\noindent \textbf{Proof.} We will use always in the proof the convention
that $\Lambda ,\Lambda _{1},\Lambda _{2}\in \mathcal{O}_{c}(X)$ with $%
\Lambda _{1}\subset \Lambda _{2}$. First we verify Condition \ref{3eq47} of
Definition \ref{3eq48}. This can easily be done as follows 
\begin{eqnarray*}
p_{\Lambda _{2},\Lambda _{1}}\circ p_{\Lambda _{2}}(\gamma _{\omega
},m_{\omega }) &=&p_{\Lambda _{2},\Lambda _{1}}(\gamma _{\omega }\cap
\Lambda _{2},m_{\omega |\gamma _{\omega }\cap \Lambda _{2}}) \\
&=&(\gamma \cap \Lambda _{1},m_{\omega |\gamma _{\omega }\cap \Lambda _{1}})
\\
&=&p_{\Lambda _{1}}(\gamma _{\omega },m_{\omega }),
\end{eqnarray*}
for $(\gamma _{\omega },m_{\omega })\in \Omega _{X}^{M}$ which is the
desired result.

\begin{figure}
\input{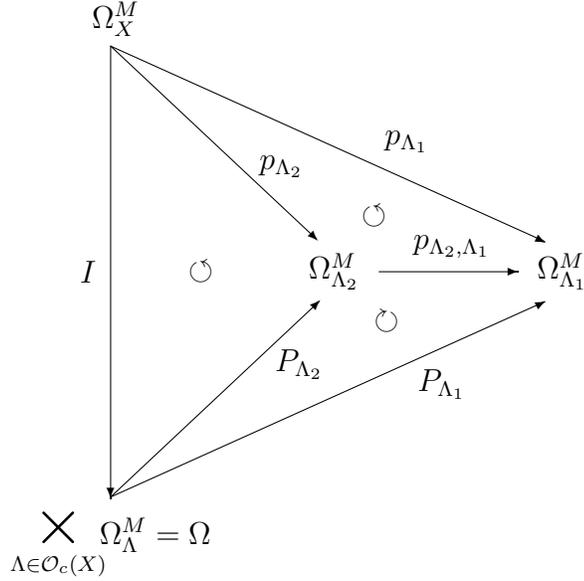}
\caption{Diagram used in constructing projective limit.\label{3eq66}}
\end{figure}

Let us now construct a version of the projective limit of the family $%
\{\Omega _{\Lambda }^{M},\Lambda \in \mathcal{O}_{c}(X)\}$, see e.g. \cite
{P67}. As $\Omega $ we take 
\[
\Omega :=\left\{ \omega \in \stackunder{\Lambda \in \mathcal{O}_{c}(X)}{{%
\mbox{\huge $\times$}}}\!\!\!\Omega _{\Lambda }^{M}|p_{\Lambda _{2},\Lambda
_{1}}((\omega )_{\Lambda _{2}})=(\omega )_{\Lambda _{1}}\right\} , 
\]
where $(\omega )_{\Lambda }$ denotes the $\Lambda $-component of $\omega $.
As projections we choose $P_{\Lambda }(\omega ):=(\omega )_{\Lambda }$ and
define the $\sigma $-algebra as $\mathcal{B}(\Omega ):=\sigma (\{P_{\Lambda
}|\Lambda \in \mathcal{O}_{c}(X)\})$, see diagram in Fig.~\ref{3eq66}.%

Now we define a (bijective) mapping $I:\Omega _{X}^{M}\to \Omega $ by 
\[
(I(\gamma _{\omega },m_{\omega }))_{\Lambda }:=p_{\Lambda }(\gamma _{\omega
},m_{\omega })\in \Omega _{\Lambda }^{M},\;(\gamma _{\omega },m_{\omega
})\in \Omega _{X}^{M},\Lambda \in \mathcal{O}_{c}(X). 
\]

We first show that $I$ is well defined, this means that 
\[
p_{\Lambda _{2},\Lambda _{1}}((I(\gamma _{\omega },m_{\omega }))_{\Lambda
_{2}})=(\gamma _{\omega }\cap \Lambda _{1},m_{\omega |\gamma _{\omega }\cap
\Lambda _{1}}). 
\]
Indeed we have 
\begin{eqnarray*}
p_{\Lambda _{2},\Lambda _{1}}((I(\gamma _{\omega },m_{\omega }))_{\Lambda
_{2}}) &=&p_{\Lambda _{2},\Lambda _{1}}\circ p_{\Lambda _{2}}(\gamma
_{\omega },m_{\omega }) \\
&=&p_{\Lambda _{1}}(\gamma _{\omega },m_{\omega }) \\
&=&(\gamma _{\omega }\cap \Lambda _{1},m_{\omega |\gamma _{\omega }\cap
\Lambda _{1}})
\end{eqnarray*}
which proves that $I$ is well defined. Let us prove in addition that $I$ is
a bijective mapping between $\Omega _{X}^{M}$ and $\Omega $.

\noindent Injectivity. Let $(\gamma _{\omega },m_{\omega })$, $(\gamma
_{\omega ^{\prime }},m_{\omega ^{\prime }})\in \Omega _{X}^{M}$ such that $%
I(\gamma _{\omega },m_{\omega })=I(\gamma _{\omega ^{\prime }},m_{\omega
^{\prime }})$, that means by definition of $I$ that $(\gamma _{\omega }\cap
\Lambda ,m_{\omega |\gamma _{\omega }\cap \Lambda })=(\gamma _{\omega
^{\prime }}\cap \Lambda ,m_{\omega ^{\prime }|\gamma _{\omega ^{\prime
}}\cap \Lambda })$ for all $\Lambda \in \mathcal{O}_{c}(X)$. Since the
manifold $X$ can be written as countable union of sets from $\mathcal{O}%
_{c}(X)$, i.e., 
\[
X=\bigcup_{n\in \Bbb{N}_{0}}\Lambda _{n},\;\Lambda _{n}\in \mathcal{O}%
_{c}(X),n\in \Bbb{N}_{0} 
\]
this implies that $(\gamma _{\omega },m_{\omega })=(\gamma _{\omega ^{\prime
}},m_{\omega ^{\prime }})$ and therefore the injectivity of $I$ is proved.

\noindent Surjectivity. Let $\omega =((\omega )_{\Lambda })_{\Lambda \in 
\mathcal{O}_{c}(X)}=((\gamma _{\omega }^{*},m_{\omega }^{*})_{\Lambda
})_{\Lambda \in \mathcal{O}_{c}(X)}\in \Omega $ be given and take a family
of pairewise disjoint subsets from $\mathcal{O}_{c}(X),$ $\{\Lambda
_{n},n\in \Bbb{N}\}$ such that 
\[
X=\bigsqcup_{n\in \Bbb{N}}\Lambda _{n}; 
\]
moreover we may assume that for any $\Lambda \in \mathcal{O}_{c}(X)$ $%
\exists m\in \Bbb{N}$ such that 
\begin{equation}
\Lambda \subset \bigsqcup_{n=1}^{m}\Lambda _{n}=:\Lambda _{n}^{m}.
\label{3eq52}
\end{equation}

Let us define an element $(\gamma ,m_{\gamma })$ from $\Omega _{X}^{M}$ as
follows: 
\[
\gamma _{\omega }:=\bigcup_{n\in \Bbb{N}}(\gamma _{\omega }^{*}\cap \Lambda
_{n}) 
\]
and 
\[
m_{\omega }:\bigcup_{n\in \Bbb{N}}(\gamma _{\omega }^{*}\cap \Lambda
_{n})\ni x\mapsto m_{x}\in M,\;x\in \gamma _{\omega }^{*}\cap \Lambda
_{n},n_{0}\in \Bbb{N}. 
\]

First we note that the assumption (\ref{3eq52}) gives 
\[
\gamma _{\omega }\cap \Lambda =\bigcup_{n=1}^{m}(\gamma _{\omega }^{*}\cap
\Lambda _{n}\cap \Lambda ),\;\Lambda \in \mathcal{O}_{c}(X). 
\]
Secondly we must prove that $(I(\gamma _{\omega },m_{\omega }))_{\Lambda
}=(\omega )_{\Lambda }$ for any $\Lambda \in \mathcal{O}_{c}(X)$. From the
definition of $I$ and $p_{\Lambda }$ we have 
\[
(I(\gamma _{\omega },m_{\omega }))_{\Lambda }:=p_{\Lambda }(\gamma _{\omega
},m_{\omega }):=(\gamma _{\omega }\cap \Lambda ,m_{\omega |\gamma _{\omega
}\cap \Lambda }) 
\]
and from the above representation for $\gamma _{\omega }\cap \Lambda $ we
obtain 
\begin{eqnarray*}
(\gamma _{\omega }\cap \Lambda ,m_{\omega |\gamma _{\omega }\cap \Lambda })
&=&\left( \bigcup_{n=1}^{m}(\gamma _{\omega }^{*}\cap \Lambda _{n}\cap
\Lambda ),m_{\omega |\gamma _{\omega }\cap \Lambda }\right) \\
&=&\bigcup_{n=1}^{m}(\gamma _{\omega }^{*}\cap \Lambda _{n}\cap \Lambda
,m_{\omega |\gamma _{\omega }^{*}\cap \Lambda _{n}\cap \Lambda }) \\
&=&\bigcup_{n=1}^{m}p_{\Lambda _{n}\cap \Lambda }(\gamma _{\omega
}^{*},m_{\omega })=\bigcup_{n=1}^{m}p_{\Lambda }(\gamma _{\omega }^{*}\cap
\Lambda _{n},m_{\omega |\gamma _{\omega }^{*}\cap \Lambda _{N}}) \\
&=&p_{\Lambda }(\gamma _{\omega }^{*}\cap \Lambda _{n}^{m},m_{\omega |\gamma
_{\omega }^{*}\cap \Lambda _{n}^{m}})=(\gamma _{\omega }^{*}\cap \Lambda
,m_{\omega |\gamma _{\omega }^{*}\cap \Lambda }) \\
&=&(\omega )_{\Lambda }
\end{eqnarray*}
which proves the surjectivity of $I$.

Taking into account Condition \ref{3eq47} of Definition \ref{3eq48} it
follows that $I$ is continuous. Hence only remained to proof that $\mathcal{B%
}(\Omega _{X}^{M})$ and $\mathcal{B}(\Omega )$ coincide. This is an
immediate consequence of the definition of the $\sigma $-algebras and the
continuity of $I$.\hfill $\blacksquare $

\subsection{Marked Poisson measure\label{3eq50}}

The underlying manifold $X$ is endowed with a non-atomic Radon measure $%
\sigma $, (cf. Section \ref{3eq46}). Let a probability measure $\tau $ be
given on the space $M$. The space $X\times M$ is endowed with the product
measure between $\sigma $ and $\tau $ denoted by $\hat{\sigma}$, i.e., $\hat{%
\sigma}:=\sigma \otimes \tau $. The measure $\hat{\sigma}^{n}$ can be
considered as a finite measure on $(\Lambda \times M)^{n}$ for any $\Lambda
\in \mathcal{O}_{c}(X)$ which induces on $\Omega _{\Lambda }^{M}(n)$ the
measure 
\[
\hat{\sigma}_{\Lambda ,n}=\frac{1}{n!}(\sigma \otimes \tau )^{n},\;n\geq 0,\;%
\hat{\sigma}_{\Lambda ,0}(\emptyset )=1. 
\]
Then we consider a measure $\lambda _{\hat{\sigma}}^{\Lambda }$ on $\Omega
_{\Lambda }^{M}$ which coincides on each $\Omega _{\Lambda }^{M}(n)$ with
the measure $\hat{\sigma}_{\Lambda ,n}$ as follows 
\[
\lambda _{\hat{\sigma}}^{\Lambda }=\sum_{n=0}^{\infty }\hat{\sigma}_{\Lambda
,n}=\sum_{n=0}^{\infty }\frac{1}{n!}(\sigma \otimes \tau )^{n}. 
\]
Measure $\lambda _{\hat{\sigma}}^{\Lambda }$ is a finite measure on $\Omega
_{\Lambda }^{M}$ and $\lambda _{\hat{\sigma}}^{\Lambda }(\Omega _{\Lambda
}^{M})=e^{\sigma (\Lambda )}$, therefore we define a probability measure $%
\mu _{\hat{\sigma}}^{\Lambda }$ on $\Omega _{\Lambda }^{M}$ setting 
\begin{equation}
\mu _{\hat{\sigma}}^{\Lambda }=e^{-\sigma (\Lambda )}\lambda _{\hat{\sigma}%
}^{\Lambda }.  \label{3eq54}
\end{equation}
The measure $\mu _{\hat{\sigma}}^{\Lambda }$ has the following property 
\[
\mu _{\hat{\sigma}}^{\Lambda }(\Omega _{\Lambda }^{M}(n))=\frac{1}{n!}\sigma
^{n}(\Lambda )e^{-\sigma (\Lambda )} 
\]
which gives the probability of the occurrence of exactly $n$ points of the
marked Poisson process (with arbitrary values of marks) inside the volume $%
\Lambda $.

\begin{figure}
\input{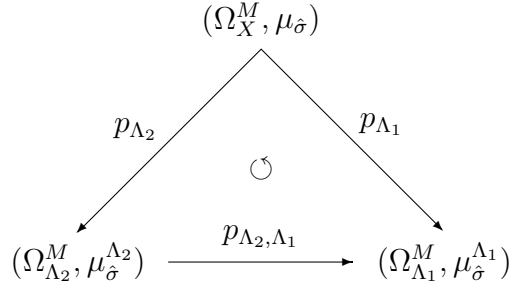}
\caption{Diagram used to prove the consistency property of the family 
$\{\mu_{\hat{\sigma}}^\Lambda | \Lambda\in{\cal O}_{c}(X)\}$.\label{3eq67}}
\end{figure}

In order to obtain the existence of a unique probability measure $\mu _{\hat{%
\sigma}}$ on $\mathcal{B}(\Omega _{X}^{M})$ such that 
\[
\mu _{\hat{\sigma}}^{\Lambda }=p_{\Lambda }^{*}\mu _{\hat{\sigma}},\;\Lambda
\in \mathcal{O}_{c}(X),
\]
one should check the consistency property of the family $\{\mu _{\hat{\sigma}%
}^{\Lambda }|\Lambda \in \mathcal{O}_{c}(X)\}$. In other words one should
verify the following equality of measures, see diagram in Fig.~\ref{3eq67}. 
\begin{equation}
\mu _{\hat{\sigma}}^{\Lambda _{2}}\circ p_{\Lambda _{2},\Lambda
_{1}}^{-1}=\mu _{\hat{\sigma}}^{\Lambda _{1}},\;\Lambda _{1},\Lambda _{2}\in 
\mathcal{O}_{c}(X),\Lambda _{1}\subset \Lambda _{2}.  \label{3eq53}
\end{equation}

It is known, see e.g.~\cite{GGV75}, that the $\sigma $-algebra $\mathcal{B}%
(\Omega _{\Lambda }^{M})$ coincides with the $\sigma $-algebra generated by
the cylinder sets from $\Omega _{\Lambda }^{M}$, $C_{B,n}^{\Lambda }$, $B\in 
\mathcal{O}_{c}(\Lambda )$, $n\in \Bbb{N}_{0}$. Here $C_{B,n}^{\Lambda }$
has the following representation 
\[
C_{B,n}^{\Lambda }:=\{\omega =(\gamma _{\omega },m_{\omega })\in \Omega
_{\Lambda }^{M}||\gamma _{\omega }\cap B|=n\}. 
\]
Hence for a given $B\in \mathcal{O}_{c}(\Lambda _{1})$, $n\in \Bbb{N}_{0}$
the pre-image under $p_{\Lambda _{2},\Lambda _{1}}$of the cylinder set $%
C_{B,n}^{\Lambda _{1}}$ from $\Omega _{\Lambda _{1}}^{M}$ is a cylinder set
from $\Omega _{\Lambda _{2}}^{M}$, i.e., 
\[
p_{\Lambda _{2},\Lambda _{1}}^{-1}(C_{B,n}^{\Lambda _{1}})=\{\omega =(\gamma
_{\omega },m_{\omega })\in \Omega _{\Lambda _{2}}^{M}||\gamma _{\omega }\cap
B|=n\}=C_{B,n}^{\Lambda _{2}}. 
\]
On the other hand it is well known, see e.g.~\cite{AKR97} and \cite{S94}
that 
\[
\mu _{\hat{\sigma}}^{\Lambda _{2}}(C_{B,n}^{\Lambda _{2}})=\frac{1}{n!}%
\sigma ^{n}(B)e^{-\sigma (B)} 
\]
which is the same as $\mu _{\hat{\sigma}}^{\Lambda _{1}}(C_{B,n}^{\Lambda
_{1}})$. Therefore the consistency property (\ref{3eq53}) is proved.

It is possible to compute in closed form the Laplace transform of the
measure $\mu _{\hat{\sigma}}$. Namely, let $f$ be a continuous function on $%
X\times M$ such that the supp$f\subset \Lambda \times M$ for some $\Lambda
\in \mathcal{O}_{c}(X)$. Let $\omega =(\gamma _{\omega },m_{\omega })$ be an
element of $\Omega _{X}^{M}$ and define the pairing between $f$ and $\omega $
by 
\[
\langle f,\omega \rangle :=\sum_{x\in \gamma _{\omega }}f(x,m_{x}). 
\]
Then we have 
\[
\int_{\Omega _{X}^{M}}e^{\langle f,\omega \rangle }d\mu _{\hat{\sigma}%
}(\omega )=\int_{\Omega _{\Lambda }^{M}}e^{\langle f,\omega \rangle }d\mu _{%
\hat{\sigma}}\circ p_{\Lambda }(\omega )=\int_{\Omega _{\Lambda
}^{M}}e^{\langle f,\omega \rangle }d\mu _{\hat{\sigma}}^{\Lambda }(\omega ). 
\]
Using (\ref{3eq54}) the last integral is equal to 
\begin{eqnarray*}
&&e^{-\sigma (\Lambda )}\sum_{n=0}^{\infty }\frac{1}{n!}\int_{(\Lambda
\times M)^{n}}\exp \left( \sum_{k=0}^{n}f(x_{k},m_{x_{k}})\right) d\hat{%
\sigma}(x_{1},m_{x_{1}})\ldots d\hat{\sigma}(x_{n},m_{x_{n}}) \\
&=&e^{-\sigma (\Lambda )}\sum_{n=0}^{\infty }\frac{1}{n!}\left(
\int_{\Lambda \times M}e^{f(x,m_{x})}d\hat{\sigma}(x,m_{x})\right) ^{n} \\
&=&\exp \left( \int_{X\times M}(e^{f(x,m_{x})}-1)d\hat{\sigma}%
(x,m_{x})\right) .
\end{eqnarray*}
That is, for any $f$ in the above conditions the following formula holds: 
\[
l_{\mu _{\hat{\sigma}}}(f)=\int_{\Omega _{X}^{M}}e^{\langle f,\omega \rangle
}d\mu _{\hat{\sigma}}(\omega )=\exp \left( \int_{X\times
M}(e^{f(x,m_{x})}-1)d\hat{\sigma}(x,m_{x})\right) . 
\]
\bigskip

\noindent Acknowledgments\medskip

We would like to thank Prof.~Dr.~R.~Minlos and Prof.~Dr.~M.~R\"{o}ckner for
helpful discussions during the preparation of this work. Financial support
of the DFG through the project AL 214/9-2, JNICT-plurianual I\&D N.~219/94,
and TMR Nr.~ERB4001GT957046 are gratefully acknowledged.

\end{document}